\documentclass[a4paper,12pt]{article}

\usepackage{amsmath}
\usepackage{latexsym}
\usepackage{epsfig}
\usepackage{amsfonts,amssymb}
\usepackage{subfigure}
\usepackage{fancybox} 

\newcommand{\reals}{\Bbb{R}}
\newcommand{\nats}{\Bbb{N}}
 
\newtheorem{theorem}{Theorem}
\newtheorem{definition}{Definition} 
 
\newtheorem{lemma}{Lemma}
\newtheorem{rem}{Remark} 
\newtheorem{prop}{Proposition} 
\newtheorem{cor}{Corollary}
\newtheorem{ass}{Assumption} 
\newtheorem{expl}{Example}
\newtheorem{claim}{Claim}

\newcommand{\N}{\mathbb{N}}
\newcommand{\ints}{\mathbb{Z}}
\newcommand{\norm}[1]{\left\vert#1\right\vert}
\newcommand{\abs}[1]{\left\vert #1 \right\vert}
\newcommand{\lab}{\label}

\newcommand{\shiftright}[1]{\null\hspace{#1}}
\newcommand{\shiftleft}[1]{\null\hspace{-#1}}
\def\rref#1{(\ref{#1})}
\newcommand\mR{\mathbb{R}}
\newcommand\mZ{\mathbb{Z}}
\newcommand{\tonio}[2][0.8\textwidth]{%
  \begin{center} %
       \begin{minipage}[c]{#1} %
         {\footnotesize\bf T:} {\footnotesize\bf #2}  %
       \end{minipage} 
  \end{center} 
}

\def\nn{\nonumber}
\def\ep{\varepsilon}
\def\dty{\displaystyle}
\newcounter{remfootnote}\addtocounter{remfootnote}{191}
\newcommand{\remfootnote}[1]{\stepcounter{remfootnote}%
                             \renewcommand{\thefootnote}{\ding{\value{remfootnote}}}%
                             \footnote{{\bf Draft comment. }#1}%
                             \renewcommand{\thefootnote}{\arabic{footnote}}}
\def\col{\mbox{ col}}
\newcommand{\putfig}[4]{
\begin{figure}[h]
\vbox{
       \epsfxsize = \textwidth 
       \divide\epsfxsize by 100
       \multiply\epsfxsize by #3
       \centerline{
       \epsffile{#1.eps} }            
      }
\caption{#4}             
\label{#2}
\end{figure}
}
\def\simulink{{\sc Simulink}$^{{\tiny {\mbox{TM}}}}$}
\def\matlab{{\sc Matlab}$^{{\tiny {\mbox{TM}}}}$}
\def\cKinfty{{\cal K}_{\infty}}

\newcommand{\noremftn}{ \renewcommand{\remfootnote}[1]{} }
\newcommand{\remappendix}{\newpage\appendix\centerline{\Large \bf Appendices for the draft}}

\usepackage{color}
\noremftn
\renewcommand{\remappendix}{\end{document}} 
\renewcommand{\tonio}[2][0.8\textwidth]{}  
\begin{document}
\title{On uniform asymptotic stability of time-varying parameterized  discrete-time cascades\thanks{This work is supported by the   Australian Research Council under the Large Grants Scheme.}
}
\author{%
\begin{minipage}{\textwidth}
\shiftright{1.3in}\begin{tabular}{c}
Dragan Ne\v{s}i\'{c}\\[2mm]
\normalsize Dept.  Electr.  Eng.,\\[-1mm]  
\normalsize The University of Melbourne, \\[-1mm] 
\normalsize Parkville, 3010, Victoria,\\[-1mm] 
\normalsize Australia.\\[-1mm] \
\end{tabular}
\begin{tabular}{c}
 Antonio Lor\'{\i}a\\[2mm]
\normalsize C.N.R.S \\[-1mm] 
\normalsize LSS - Sup\'elec,\\[-1mm] 
\normalsize Plateau de Moulon, \\[-1mm] 
\normalsize 91192, Gif sur Yvette, \\[-1mm] 
\normalsize France.
\end{tabular} 
\end{minipage}}

\date{  
 }
\maketitle

\begin{abstract}
Recently, a framework for controller design of sampled-data nonlinear systems via their approximate discrete-time models has been proposed in the literature. In this paper we develop novel tools that can be used within this framework and that are very useful for tracking problems. In particular, results for stability analysis of parameterized time-varying discrete-time cascaded systems are given. This class of models arises naturally when one uses an approximate discrete-time model to design a stabilizing or tracking controller for a sampled-data plant. While some of our results parallel their continuous-time counterparts, the stability properties that are considered, the conditions that are imposed and the the proof techniques that are used are tailored for approximate discrete-time systems and are technically different from those in the continuous-time context. We illustrate the utility of our results in the case study of the tracking control of a mobile robot. This application is fairly illustrative of the technical differences and obstacles encountered in the analysis of discrete-time parameterized systems. 
\end{abstract}

\section{Introduction}


The prevalence of digitally controlled systems and the fact that the nonlinearities in the plant model can often not be neglected, strongly motivate the area of nonlinear sampled-data systems. A typical nonlinear sampled-data system consists of a nonlinear continuous-time plant and a nonlinear discrete-time controller that are interconnected via the analog-to-digital (A-D) and digital-to-analog (D-A) converters. Despite the importance of this class of systems, few systematic tools for nonlinear sampled-data controller design are available in the literature. {Instead, in control practice one typically follows the commonly accepted hypothesis  that if a continuous-time controller designed for a continuous-time plant, is implemented at a sufficiently fast sampling rate, the sampled-data system should ``behave well''. While this intuition is correct in general, the required sampling may be too fast to be implemented in practice because of the available hardware limitations.}

{The stumbling block that naturally arises in the} 
formal analysis and controller design for nonlinear sampled-data systems is the fact that the model of the system is rather complex (hybrid, nonlinear, periodically time-varying) and very hard to deal with directly. Consequently, there are several different methods that can be used when designing the controller. 

{One method},
which is sometimes referred to as {\em the emulation method}, 
{consists of designing a continuous-time controller} based on the continuous-time plant model and then discretizing the controller for digital implementation. {We stress } that the emulation method ignores {the} sampling during the controller design step. Results on emulation for nonlinear sampled-data systems can be found in \cite{dina-nesic-teel} and references defined therein. 
{Another method}, {which for ease of reference we will call here, exact (respectively approximate) discrete-time design (DTD),}  {consists of
obtaining the exact (respectively, approximate) discrete-time model of the plant} and then designing the controller based on the discrete-time plant model. In this method,   the inter-sample behaviour is ignored during the controller design. 

{Since the emulation method does not take the sampling into account in the controller design step, it is reasonable to expect that DTD method may produce better results than the emulation. Unfortunately, if we want to implement the second method}, even if the nonlinear continuous-time plant model is known, in general we will be unable to obtain the exact discrete-time model of the plant since this involves solving a nonlinear differential equation analytically over one sampling interval. Instead, we need to exploit an approximate discrete-time model that is obtained using some numerical integration scheme, such as Runge-Kutta (see e.g. \cite{gruene-kloeden,stuart}). Early results that use approximate discrete-time plant models for controller design can be found in \cite{dochain,goodwin}. See also the more recent works of \cite{monaco1,monaco2,mareels}. 

{The main pitfall of the approximate DTD method is that if one is not careful with the choice of the approximate model and the design of the controller, it is possible that a controller asymptotically stabilizes the approximate plant model but {\em not} the exact model. It is noteworthy that this fact concerns even linear systems as we illustrate next.
\begin{expl}\cite{NESTEEPK}
Consider the system
\begin{equation}\label{doubleint}
\left\{ \, \begin{array}{lcl}
\dot x_1 &=& x_2\\
\dot x_2 &=& u
\end{array}\right.
\end{equation}
whose exact and Euler models are respectively
\begin{eqnarray}
\mbox{Euler:} & \left\{ \
\begin{array}{lcl}
x_1(k+1) & = & x_1(k) + Tx_2(k)\\
x_2(k+1) & = & x_2(k) + Tu(k)
\end{array}\right.   \\[2mm]
\mbox{Exact:} & \left\{ \  
\begin{array}{lcl}
x_1(k+1) & = & x_1(k) + Tx_2(k) + 0.5T^2u(k)\\
x_2(k+1) & = & x_2(k) + Tu(k)\,.
\end{array}\right. 
\end{eqnarray}
If we define the controller $u_T(x)= -(1/T)[ x_1 + 2x_2 ]$  for the Euler model, we get that one eigen-value of the closed loop system with the exact-model is on the unitary circle for all values of $T>0$. Hence, the system cannot be asymptotically stable for any $T>0$. On the other hand, the eigenvalues of the approximate model are $\pm \sqrt{(1-T)}$ and one can actually show that there exists $b>0$ such that $\norm{x(k)} \leq b\norm{x(0)}\mbox{e}^{-0.5kT}$ for all $x(0)\in\mR^2$ and all $T\in (0,0.5)$. 
\end{expl}
}

Motivated by this fact, a framework for nonlinear sampled-data controller design via approximate discrete-time models has been proposed in \cite{NESTEEPK,NESTEESON,nesic-teel,NESLAI,nesic-angeli}. These results are very similar in spirit to results from the numerical analysis literature (see e.g. \cite{stuart,lars-book}) that applies to continuous-time control systems.  {In \cite{NESTEEPK}} checkable conditions on the continuous-time plant model, the controller and the approximate discrete-time model are presented which guarantee that if the controller stabilizes the approximate model, it would also stabilize {the exact discrete-time model. 
 Furthermore, in \cite{NESTEESON} it was shown that stability of the exact discrete-time model under mild conditions guarantees also stability of the sampled-data system. Hence, the results of \cite{NESTEEPK,NESTEESON} provide a framework for controller design of sampled-data nonlinear systems via their approximate discrete-time models. }


The above mentioned results are {primarily targeted at establishing conditions for stability that is, they} are {\em prescriptive} and {\em non-constructive}. {In other words,} they provide a framework for controller design {but} without explicit recipes for controller design.
 A range of different constructive methods for controller design within the above given framework has been reported in the literature: backstepping via the Euler model of strict feedback systems \cite{nesic-teel}; Lyapunov methods based on changes of supply rates for input-to-state stable systems \cite{iss-nesic-teel,dina-nesic,sg-laila-nesic}; and optimization based stabilization \cite{gruene-nesic}. 

Other constructive design methods for non-parameterized systems can be found, for instance in \cite{MADMONNOR,mareels,LINBYR,SIMNIJTSI}, in  the survey \cite{monaco2} and in the references listed therein.

{\subsubsection*{Contributions of this paper}
The results that we present in this paper contribute to what we may call {\em cascades-based} control. Roughly speaking,  this approach aims at  designing controllers in cases when the closed loop system has a cascaded structure. Moreover, the closed loop dynamics shall verify certain structural conditions imposed either on the functions that define the closed loop dynamics or indirectly, in terms of properties of Lyapunov-like functions. While there is a wide number of such results establishing different forms of asymptotic stability for continuous-time systems (both, autonomous and non-autonomous), there are only a few results for discrete-time systems. In particular, we cite \cite{JIAWAN} where sufficient conditions for stability of cascades that use the  input-to-state stability (ISS) property are presented. }

{Yet, the study of stability of cascaded systems in the discrete-time context has a double motivation: Firstly, it obviously inherits the motivations from the continuous-time case: there is a wide range of applications in control design, of stability results for cascades.
See \cite{MURAT-CASC,ERJENSTHESIS,CASCIPN,CASCAUT} and references therein for a {large number} of results and applications in this area; }
from a control theory viewpoint, the motivation for this research originated probably in geometric nonlinear control where it was shown that many systems can be transformed into a cascade via a local change of coordinates (see, for example, \cite[Lemma 1.6.1]{ISI}). {Secondly, in view of the previous discussions it is highly desirable to establish conditions tailored specifically for approximate discrete-time systems, under which one can rely on the DTD control design method.}

This paper {is aimed in that direction. We } present results on uniform asymptotic practical stability of parameterized (in the sampling period) discrete-time cascaded  systems. Thereby, contributing to the framework established in  \cite{NESTEEPK,nesic-teel,NESTEESON} and other above-cited references for sampled-data systems. Although 
we have been inspired by 
 continuous-time results on stability of cascades (in particular by those in \cite{MURAT-CASC,CASCAUT}), the results presented here are not a simple translation of their counterparts in continuous-time. Indeed, the properties we consider, the conditions we impose and the proofs we establish here are notably different from continuous-time ones and they are tailored specifically for discrete-time parameterized systems that arise within the above mentioned framework. 

{ Furthermore, advantages of DTD control design are illustrated with the sampled-data tracking control problem of the unicycle benchmark problem where we show an improvement in performance with respect to continuous-time based designs. In particular, we will illustrate how to apply our main results on stability of cascades. While the benchmark problem is the same as in continuous time, this case-study is fairly representative of the large differences between discrete-time and continuous-time based designs as well as the tools involved in the  study of the proper stability properties for parameterized systems. Indeed, the proofs of the results for the unicycle are original and also very different from their continuous counter-parts. Thus, through this case-study we illustrate how using our results we can obtain new control algorithms that can be regarded as a continuous-time controller redesign  that outperforms the emulated controllers.}

The rest of the paper is organized as follows. In Section \ref{sec:next} we provide mathematical preliminaries with some definitions and basic results for parameterized systems. Our main results are presented in Section \ref{sec:main} in the form of two different theorems that are ``Lyapunov-based'' and ``trajectory based'' respectively. 
In Section \ref{sec:appl} we present the unicycle case study and the proofs of the main theorems are presented in Section \ref{sec:proofs}. Some concluding remarks are presented in Section \ref{sec:concl}. Further technical proofs and other auxiliary results are included in the appendix.

\noindent {\bf Notation}. A function $\alpha : \reals_{\geq 0} \rightarrow \reals_{\geq 0}$ is said to be of class ${\cal K}$ ($\alpha \in {\cal K}$), if it is continuous, strictly increasing and zero at zero; $\alpha \in {\cal K}_\infty$ if, in addition, it is unbounded.  A function $\beta : \reals_{\geq 0} \times \reals_{\geq 0} \rightarrow \reals_{\geq 0}$ is of class ${\cal K}{\cal L}$ if for all $t > 0$, $\beta(\cdot,t)\in {\cal K}$, for all $s > 0$, $\beta(s,\cdot)$ is decreasing to zero. A function $\gamma: \reals_{\geq 0} \rightarrow \reals_{\geq 0}$ is said to be of class ${\cal N}$ if $\gamma(\cdot)$ is continuous and nondecreasing. We denote by $\norm{\cdot}$ the Euclidean norm of vectors. We denote by $\reals$ and $\nats$ the sets of the real and natural numbers respectively. For an arbitrary $r \in \reals$ we use the notation $\lfloor r \rfloor:= \dty\max_{z \in \ints, z \leq r} z$. Given strictly positive real numbers $L,T$ we use the following notation:
\begin{equation}
\label{ell}
\ell_{L,T}:= \left\lfloor \frac{L}{T} \right\rfloor \ .
\end{equation}


\section{{Parameterized discrete-time systems}}
\label{sec:next}

In this section we present a result that provides a framework for controller design for sampled-data systems via their approximate discrete-time models. This result generalizes \cite[Theorem 1]{NESTEEPK} in that it is applicable to time-varying plants and hence more suitable for tracking problems. The proof of this result follows the same steps as that of Theorem 1 in the above paper and hence it is omitted. Our main results that are presented in Section 3 facilitate sampled-data controller design within the framework of Theorem 1 in this section. 

We also introduce here the precise definitions of stability that we pursue and other preliminary important results.

\subsection{Controller design framework}

Consider the class of systems
\begin{eqnarray} \lab{orig}
	\dot{x}(t) & = & f(t,x(t),u(t)) \nn \\
     	y(t) & = & h(x(t))
\end{eqnarray}
where $x \in {\mR}^{n_x}$ and $u \in {\mR}^{m}$ are respectively the state and control input. 
We assume that for any given $x_\circ$ and $u(\cdot)$ the differential equation in (\ref{orig}) has a unique solution defined on its maximal interval of existence $[0,t_{\max})$. This may be guaranteed, for instance, by requiring $f$ in (\ref{orig}) to be locally Lipschitz, uniformly in $t$.
The control is taken to be a piecewise constant signal $u(t) = u(kT) =:u(k) , \ \forall t \in [kT,(k+1)T)$, $k \in \N$, where $T > 0$ is the sampling period. Assume that some combination (output) or all of the states ($x(k) := x(kT)$) are available at sampling instant $kT, k \in \N$. The exact discrete-time model for the plant (\ref{orig}), which describes the plant behavior at sampling instants $kT$, is obtained by integrating the 
initial value problem 
\begin{align} \lab{dts}
	\dot{x}(t)=f(t,x(t),u(k)) \ ,  
\end{align} 
with given $u(k)$ and $x_0=x(k)$, over the sampling interval $[kT,(k+1)T]$. 

If we denote by $x(t)$ the solution of the initial value problem (\ref{dts}) at time $t$ with given $x_0=x(k)$ and $u(k)$, then the exact discrete-time model of (\ref{orig}) can be written as:
\begin{equation} \lab{exact} 
\begin{split}
	x(k+1) &= x(k)+\int_{kT}^{(k+1)T}f(\tau,x(\tau),u(k))d\tau
	=: F_T^{e}(x(k),u(k)) \ .
\end{split}
\end{equation}
We emphasize that $F_T^{e}$ is not known in most cases. Indeed, in order to compute $F_T^e$ we have to solve the initial value problem (\ref{dts}) analytically and this is usually impossible since $f$ in (\ref{orig}) is nonlinear. Hence, we will use an approximate discrete-time model of the plant to design a discrete-time controller for the original plant (\ref{orig}). 

Different approximate discrete-time models can be obtained using different methods, such as a classical Runge-Kutta numerical integration scheme for the initial value problem (\ref{dts}) --see e.g. \cite{gruene-kloeden,monaco1,monaco2,stuart,ferretti}. The approximate discrete-time model can be written as
\begin{equation}\lab{approx}
	x(k+1) = F_T^a(k,x(k),u(k)) \ .
\end{equation}
For instance, if $f$ is locally Lipschitz in $t$ and $x$, the Euler approximate model can be defined as $x(k+1)=x(k)+Tf(kT,x(k),u(k))$ and it can be shown to be an $O(T^2)$ approximation of the exact discrete-time model. On the other hand, if $f$ is measurable in $t$, then a modified ``Euler'' model that is $O(T^2)$ approximation of the exact model is given by $x(k+1)=x(k)+\int_{kT}^{(k+1)T}f(\tau,x(k),u(k)) d\tau$. 

In our work, the sampling period $T$ is assumed to be a design parameter which can be {\em arbitrarily assigned}. Since we are dealing with a family of approximate discrete-time models $F_T^a$, {\em parameterized} by $T$, in order to achieve a certain objective we need in general to obtain a family of controllers, parameterized by $T$. Thus, we consider a family of dynamic feedback controllers
\begin{equation} \lab{ctrl}
\begin{array}{rcl}
  z(k+1) & = & G_T(k,x(k),z(k))  \\
	u(k) & = & u_T(k,x(k),z(k)) \ ,
\end{array}
\end{equation}
where $z \in {\mR}^{n_z}$. We denote the right hand sides of the closed-loop systems (\ref{approx}), (\ref{ctrl}) and (\ref{exact}), (\ref{ctrl}) respectively as 
$$
\widetilde{F}_T^a(k,\tilde{x}):= \left(
\begin{array}{c}
F_T^a(k,x,u_T(k,x,z)) \\
G_T(k,x,z)
\end{array}
\right) ,
\qquad
\widetilde{F}_T^e(k,\tilde{x}):=
\left(
\begin{array}{c}
F_T^e(k,x,u_T(k,x,z)) \\
G_T(k,x,z)
\end{array}
\right)
$$
where $\tilde{x}=(x^T \ z^T)^T$. 

We emphasize again that if the controller (\ref{ctrl}) stabilizes the approximate model (\ref{approx}) for all small $T$, this does not guarantee that the same controller would approximately stabilize the exact model (\ref{exact}) for all small $T$. In addition to Example \ref{doubleint}, several inherently different examples that illustrate this phenomenon can be found in \cite{NESTEEPK,nesic-teel,gruene-nesic,sg-laila-nesic}. The following result that was proved in \cite{NESTEEPK} for time-invariant systems, gives sufficient conditions for stabilization of (\ref{exact}) via controllers (\ref{ctrl}) that are designed using (\ref{approx}). The result states that one can conclude semiglobal uniform asymptotic stability (SP-UAS) of the exact discrete-time closed loop system from the same property for the approximate discrete-time closed loop and an appropriate consistency property between the exact and approximate models (see Definition \ref{def:1}). We emphasize that under very weak conditions this guarantees that the inter-sample behaviour of the sample-data closed-loop system would also be bounded (see \cite{NESTEESON}). 

The following result generalizes \cite[Theorem 1]{NESTEEPK} to time-varying systems. The proof is omitted since it follows from that of Theorem 1 in the above paper with straightforward changes.
\begin{theorem} 
\lab{th-motivate}
Suppose that there exists $\beta \in {\cal KL}$ such that for any strictly positive numbers $(L,\eta,\Delta,\delta)$ there exist $\alpha: \mR_{\geq 0} \times \mR_{\geq 0} \rightarrow \mR_{\geq 0} \cup \infty$ and $T^*>0$ such that for all $k_\circ \geq 0$, $\abs{\tilde{x}(k_\circ)} \leq \Delta$ and $T \in (0,T^*)$ we have:
\begin{itemize}
\item {\bf SP-UAS of approximate:} \\
The solutions of (\ref{approx}), (\ref{ctrl}) satisfy:
\begin{eqnarray} \label{approx-bnd}
	\abs{\phi_{T}^a(k,k_\circ,\tilde{x}(k_\circ))} & \leq & \beta(\abs{\tilde{x}(k_\circ)},T(k-k_\circ))+\delta , \qquad k \geq k_\circ \geq 0 \ .
\end{eqnarray}
\item {\bf multiple-step consistency between $\widetilde{F}_T^a$ and $\widetilde{F}_T^e$:} \\
$\max\{\tilde{x}_1,\tilde{x}_2 \} \leq \Delta$, $\abs{\tilde{x}_1-\tilde{x}_2} \leq c$ and $k \geq 0$ imply
\begin{equation}
\abs{\widetilde{F}_T^e(k,\tilde{x}_1) - \widetilde{F}_T^a(k,\tilde{x}_2)} \leq \alpha(c,T)
\end{equation}
where $\alpha$ satisfies
\begin{equation}
\alpha^k(0,T) \leq \eta , \ \qquad \forall k \in \left[0,\ell_{L,T} \right]
\end{equation}
\end{itemize}
Then, for any strictly positive real numbers $(\widetilde{\Delta},\widetilde{\delta})$ there exists $\widetilde{T}>0$ such that for all $k_\circ \geq 0$, $\abs{\tilde{x}(k_\circ)} \leq  \widetilde{\Delta}$ and $T \in (0,\widetilde{T})$ the solutions of (\ref{exact}), (\ref{ctrl}) satisfy:
\begin{itemize}
\item {\bf SP-UAS of exact:} \\
\begin{equation}
\abs{\phi_{T}^e(k,k_\circ,\tilde{x}(k_\circ))} \leq \beta(\abs{\tilde{x}(k_\circ)},T(k-k_\circ))
	+\widetilde{\delta}, \qquad \forall k \geq k_\circ \geq 0 \ . 
\end{equation}
{\scriptsize \null \hfill \null$\square$}
\end{itemize}
	
\end{theorem}
Checking SP-UAS for the approximate model (see the equation \rref{approx-bnd}) is very hard in general and in this paper we provide new results that facilitate checking this property in the case when the approximate model has a cascaded structure.
\begin{rem}
\label{rem:consistency}
The consistency condition in Theorem \ref{th-motivate} is checkable although $F_T^e$ is not known in general. For instance, it was shown that this condition holds if for any $\Delta>0$ there exist $K,T^*>0$ and $\rho \in {\cal K}$ such that $\max\{\abs{\tilde{x}_1},\abs{\tilde{x}_2} \} \leq \Delta$, $k \geq 0$ and $T \in (0,T^*)$ imply
$$
\abs{\widetilde{F}_T^e(k,\tilde{x}_1)-\widetilde{F}_T^a(k,\tilde{x}_2)} \leq (1+KT) \abs{\tilde{x}_1-\tilde{x}_2}+ T\rho(T) \ .
$$
This condition is readily checkable and it is commonly used in numerical analysis literature \cite{stuart}. For instance, it holds if the following two inequalities hold in a semiglobal practical sense:
\begin{eqnarray}
\abs{\widetilde{F}_T^a(k,\tilde{x}_1)-\widetilde{F}_T^a(k,\tilde{x}_2)} & \leq &  (1+KT) \abs{\tilde{x}_1-\tilde{x}_2} \label{cond-lip} \\
\abs{\widetilde{F}_T^e(k,\tilde{x})-\widetilde{F}_T^a(k,\tilde{x})} & \leq & T\rho(T) \ . \label{cond-cosist}
\end{eqnarray}
The condition (\ref{cond-lip}) is a particular local Lipschitz property of $\widetilde{F}_T^a$, which is usually not hard to check (see Remark \ref{rem:proofUGC:cs}). The condition (\ref{cond-cosist}) is ensured by choosing an appropriate consistent approximate model $F_T^a$ for controller design. A range of such approximate models can be found in standard books on numerical analysis \cite{stuart} (if $f$ is independent of $t$), in \cite{ferretti} (if $f$ is differentiable in $t$) and in \cite{gruene-kloeden} (if $f$ is measurable in $t$ in the Lebesgue sense).
\end{rem}

Theorem \ref{th-motivate} and the discussion in the Introduction motivate the following stability property. Note that below we use the ``max'' instead of the ``$+$'' characterization of the practical stability that was used in Theorem \ref{th-motivate} since the ``max'' formulation is easier to use in our proofs, but the two characterizations are qualitatively the same.
\begin{definition}[SP-UAS]
\label{def:1}
The family of the parameterized time-varying systems
\begin{equation}
\label{system-whole}
y(k+1)=F_T(k,y(k))
\end{equation}
is semiglobally practically uniformly asymptotically stable SP-UAS (resp. uniformly globally asymptotically stable UGAS) if there exists $\beta \in {\cal KL}$ such that for any pair of strictly positive real numbers $(\Delta,\nu)$ there exists $T^*>0$ (resp. there exists $T^*>0$) such that for all $k_\circ \geq 0$, $y(k_\circ)=y_\circ$ with $\abs{y_\circ} \leq \Delta$, $T \in (0,T^*)$ ($y(k_\circ)=y_\circ$ with $y_\circ \in \reals^n$, $T \in (0,T^*)$) the following holds:
\begin{eqnarray} 
\label{bound}
&& \shiftleft{0.8cm} \dty\abs{\phi^y_T(k,k_\circ,y_\circ) } \leq \max\{ \beta(\abs{y_\circ},(k-k_\circ)T), \nu \} \qquad (\mbox{resp.} \ \abs{\phi^y_T(k,k_\circ,y_\circ)} \leq \beta(\abs{y_\circ},(k-k_\circ)T) \ ) \nonumber  \\  &&
\end{eqnarray}
for all $k \geq k_\circ$. \null \hfill \null $\square$
\end{definition}
Notice that the convergence of solutions in (\ref{bound}) is not uniform in the sampling period since the solutions are allowed to converge slower as $T$ increases. However, the bound \rref{bound}  has  uniform overshoots. The motivation to consider this precise type of asymptotic stability is twofold. First, this is the same property that the approximate model should satisfy in the statement of Theorem \ref{th-motivate}. Second, under the conditions of this theorem, the exact model would also satisfy the stability property of Definition \ref{def:1} and by further applying results of \cite{NESTEESON} we can conclude that the sampled-data systems is SP-UAS.

\subsection{Problem setting and other definitions}

This paper focusses on the problem of establishing sufficient and necessary conditions for SP-UAS of cascaded parameterized systems of the general form 
\begin{eqnarray}
x(k+1) & = & f_T(k,x(k),z(k)) \lab{sys1} \\
z(k+1) & = & g_T(k,z(k))  \ , \lab{sys2}
\end{eqnarray} 
where $x\in\mR^{n_x}$, $z\in\mR^{n_z}$ and $T$ is the sampling period. For the cascade (\ref{sys1}), (\ref{sys2}) we use the notation $\xi:=[x^T \ z^T]^T$ to denote the state of the overall system. These results are used within the framework provided by Theorem \ref{th-motivate} to check appropriate stability properties of the approximate discrete-time model.

To state our main results we also consider the system
\begin{eqnarray}
x(k+1) & = & f_T(k,x(k),0) \lab{sys1b} \ ,
\end{eqnarray}
and often, we will regard $z$ in the system (\ref{sys1}) as an exogenous input that is not necessarily generated by the subsystem (\ref{sys2}). In that case, we refer to the subsystem (\ref{sys1}) as the {\em system with input $z$}. The solution of the system (\ref{sys1}) with input $z$ at time $k$ that starts at initial time instant $k_\circ$ from the initial state $x(k_\circ)=x_\circ$ and under the action of the input sequence $\omega^z_{[k_\circ,k)}:=\{z(k_\circ),\ldots z(k-1)\}$ is denoted as $\phi^x_T(k,k_\circ,x_\circ,\omega^z_{[k_\circ,k)})$. We also use $\omega^z:=\omega^z_{[k_\circ,\infty)}$. Note that the solution of the system (\ref{sys1b}) is the same as the solution for system (\ref{sys1}) with input $z$ when $z(j) \equiv 0 , \forall j \in [k_\circ,k]$ and hence for solutions of (\ref{sys1b}) we use the notation $\phi^x_T(k,k_\circ,x_\circ,0)$. Similarly we use the notation $\phi_T^\xi(k,k_\circ,\xi_\circ)$, $\phi_T^x(k,k_\circ,\xi_\circ)$ and $\phi_T^z(k,k_\circ,z_\circ)$ to denote solutions of the overall system (\ref{sys1}), (\ref{sys2}) and its $x$ and $z$ components respectively.

{Roughly speaking, our main results establish SP-UAS of the cascade under the condition that both subsystems, \rref{sys2} and \rref{sys1b} are SP-UAS and that the solutions of the cascade be uniformly bounded in the following specific manner. }
\begin{definition}
\label{def:2}
The system (\ref{system-whole}) is uniformly semiglobally bounded USB (uniformly globally bounded UGB), if there exist $\kappa\in {\cal K}_\infty$ and $c$, such that for any $\Delta>0$ there exists $T^*> 0$ (there exists $T^*>0$) such that $k_\circ \geq 0$, $y(k_\circ)=y_\circ$ with $\abs{y_\circ} \leq \Delta$ and $T \in (0,T^*)$ ($y_\circ \in \reals^n$ and $T\in (0,T^*)$) implies
\begin{equation}
\label{eq:ugb}
\abs{\phi^y_T(k,k_\circ,y_\circ)} \leq \kappa(\abs{y_\circ}) + c \ ,
\end{equation}
for all $k\geq k_\circ$. \null \hfill \null $\square$
\end{definition}
We  also use the following property {which is  very related to SP-UAS.} (cf. remark \ref{rem:converse}).
\begin{definition} 
\label{def:3}
The family of the parameterized time-varying systems (\ref{system-whole}) is Lyapunov SP-UAS if there exist  $\alpha_1$, $\alpha_2$ $\in {\cal K}_\infty$, $\alpha_3\in {\cal K}$, $L \in {\cal N}$ and for each pair  $(\Delta,\nu) > 0$ there exists $T^*>0$ and for each  $T\in (0,T^*)$ a continuous function $V_T : \reals_{\geq 0} \times \reals^{n} \to \reals_{\geq 0}$ such that for all $\abs{y}\leq \Delta$, all $k\geq 0$ and all $T\in (0,T^*)$ we have that
\begin{eqnarray}
\label{newa1a:VT:1}
&\displaystyle \alpha_1(\abs{y}) \leq V_T(k,y) \leq \alpha_2(\abs{y}) &\\
\label{newa1a:VT:2}
&\displaystyle  V_T(k+1, F_{T}(k,y))  - V_T(k, y )  \leq -  T(\alpha_3(\abs{y}) + \nu)\,& 
\end{eqnarray}
and
\begin{equation}
\label{newa1a:VT:3}
  \abs{\,V_T(k,r)  - V_T(k,s)\,}  \leq L\left(\max\{\abs{r},\abs{s} \}\right) \abs{r-s}\,.
\end{equation}
for all $\max\{\abs{r},\abs{s}\} \leq \Delta$, $T \in (0,T^*)$ and $k \geq 0$. The system is Lyapunov UGAS if there exists $T^*>0$ such that the above conditions hold for all $x,r,s \in \reals^n$, $k \geq 0$ and with $\nu=0$. \null \hfill \null $\square$
\end{definition}
The Lipschitz condition \rref{newa1a:VT:3} holds in particular if $V_T(k,\cdot)$ is continuously differentiable and the derivative of $V_T(k,\cdot)$ is bounded on compact sets, uniformly in small $T$ and $k\geq 0$. 
\begin{rem} 
\label{rem:converse}
We note that the properties in Definitions \ref{def:2} and \ref{def:3} are very related. In particular, it was shown in \cite{NESTEEPK} that (\ref{newa1a:VT:1}), (\ref{newa1a:VT:2}) are equivalent to the bound (\ref{bound}). However, the converse Lyapunov theorem in \cite{NESTEEPK} does not produce a Lyapunov function satisfying the condition (\ref{newa1a:VT:3}) and we believe that constructing such converse Lyapunov functions is an open problem in the literature. In particular, it would be important to provide conditions under which one can construct Lyapunov functions from Definition \ref{def:3} for families of systems (\ref{system-whole}) with a discontinuous right hand side. Such converse theorems for non-parameterized systems can be found in \cite{KELPHD} but we are not aware of similar results for parameterized systems.
\end{rem}
The following property is similar to continuity of solutions of differential equations satisfying the local Lipschitz condition if we think of $t:=kT$ as ``continuous time''. This will be crucial in establishing a trajectory-based proof of our main result.

\begin{definition} 
\label{def:4}
The solutions of the system (\ref{sys1}) with input $z$ are uniformly semiglobally continuous USC (uniformly globally continuous UGC) if for any $\Delta$ there exists $T^*>0$ (there exists $T^*>0$) such that for any $\eta \in (0,\Delta)$, $\epsilon>0$ and $L>0$ there exists $\mu>0$ such that for all $T \in (0,T^*)$, all $z(\cdot)$ with $\norm{\omega^z} \leq \mu$, $k_\circ \geq 0$ and all $x(k_\circ)=x_\circ$ with $\abs{x_\circ} \leq \eta$ we have that:
\begin{equation}
\label{gato}
\abs{\phi_T^x(k,k_\circ,x_\circ,\omega^z_{[k_\circ,k)}) - \phi_T^x(k,k_\circ,x_\circ,0)} \leq \epsilon \ , \ \ \forall k \in [k_\circ,k_\circ+\ell_{L,T}] \ .
\end{equation}
\null \hfill \null $\square$
\end{definition}
The last technical tool that we introduce is a set of checkable conditions for USC inspired by the literature on numerical methods. The proof can be found in the Appendix.
\begin{prop}
\lab{lem3}
Suppose that the system (\ref{sys1b}) is SP-UAS (UGAS) and for any pair of strictly positive $(\Delta_1,\Delta_2)$ there exist strictly positive $K,T^*$ such that: 
For all $k\geq 0$, $\max\{ \abs{x_1}, \abs{x_2} \} \leq \Delta_1$, $\max\{ \abs{z_1}, \abs{z_2} \} \leq \Delta_2$ and $T \in (0,T^*)$ we have:
\begin{eqnarray}
\abs{f_T(k,x_1,z_1)-f_T(k,x_2,z_1)} & \leq & (1+KT)\abs{x_1-x_2} \label{cont1} \\
\abs{f_T(k,x_1,z_1)-f_T(k,x_1,z_2)} & \leq & KT \abs{z_1-z_2} \label{cont2} \,.
\end{eqnarray}
Then, the system (\ref{sys1}) with input $z$ is USC (UGC).
\null \hfill \null $\square$
\end{prop}
Note the slight difference between conditions in Proposition \ref{lem3} and conditions presented in Remark \ref{rem:consistency}. 

\section{SP-UAS of time-varying cascades}
\label{sec:main}

We present next two main results that deduce stability of the cascade \rref{sys1}, \rref{sys2} from stability of two lower dimensional auxiliary systems \rref{sys1b} and \rref{sys2} and the boundedness of trajectories of the overall system \rref{sys1}, \rref{sys2}. These results are fundamental in that they are necessary and sufficient for stability of the cascade under appropriate conditions. Hence, the first condition of Theorem \ref{th-motivate} can be checked via by three separate conditions that are usually easier to verify. Our two main results are similar in spirit but they are derived under slightly different conditions. 

Theorem \ref{th1} states that under the USC assumption, we have that USB (UGB) of the system (\ref{sys1}), (\ref{sys2}) and SP-UAS (UGAS) of the auxiliary systems (\ref{sys1b}) and (\ref{sys2}) is equivalent to SP-UAS (UGAS) of the overall cascade (\ref{sys1b}) and (\ref{sys2}). The proof of this theorem is inspired by the trajectory based proof in \cite[Theorem 1]{SON89} that does not appeal to converse Lyapunov theorems (see Section \ref{sec:proofs}). Theorem \ref{th2} uses Assumption \ref{ass2} (see below) instead of the USC assumption and it shows that USB (UGB) of the system (\ref{sys1}), (\ref{sys2}), Lyapunov SP-UAS (UGAS) of the auxiliary system (\ref{sys1b}) and SP-UAS (UGAS) of the auxiliary system (\ref{sys2}) implies SP-UAS (UGAS) of the overall cascade (\ref{sys1b}) and (\ref{sys2}). We stress  that in Theorem \ref{th2} we use a different assumption from USC which is typically weaker than the latter. Moreover, we use Lyapunov SP-UAS which is typically stronger than SP-UAS.  In particular, in Theorem \ref{th2} we may have that the right hand side of (\ref{sys1b}) is discontinuous, which is in general excluded from Theorem \ref{th1} because of the USC condition.  Note also that since appropriate converse Lyapunov theorem does not exist for parameterized systems (see Remark \ref{rem:converse}), we only state sufficiency results in Theorem \ref{th2}. Finally, we note that the proof of Theorem \ref{th2} is inspired by the proof of \cite[Lemma 2]{CASCAUT} (see Section \ref{sec:proofs}).

\begin{theorem} 
\lab{th1}
Suppose that the solutions of the system (\ref{sys1}) with input $z$ {are} USC (UGC). Then, the system (\ref{sys1}), (\ref{sys2}) is SP-UAS (UGAS) if and only if the following conditions hold:
\begin{enumerate}
\item The system (\ref{sys1b}) is SP-UAS (UGAS);
\item The system (\ref{sys2}) is SP-UAS (UGAS);
\item The system (\ref{sys1}), (\ref{sys2}) satisfies the property USB (property UGB).
\end{enumerate}
\null \hfill \null $\square$
\end{theorem}
In order to state our second main result we need the following assumption.
\begin{ass}
\label{ass2}
There exist $\gamma_2 \in {\cal N}$, $\gamma_1,\gamma_3 \in {\cal K}_\infty$ and $T^*>0$ such that for all $\xi \in \reals^n$, $k \geq 0$ and $T \in (0,T^*)$ we have:
\begin{eqnarray}
\abs{f_T(k,x,z)} & \leq & \gamma_1(\abs{\xi}) \nonumber \\
\abs{f_T(k,x,z)-f_T(k,x,0)} & \leq & T \gamma_2(\abs{x}) \gamma_3(\abs{z}) \ . \nonumber 
\end{eqnarray}
\null \hfill \null $\square$
\end{ass}

\begin{theorem} 
\lab{th2}
Suppose that $f_T$ of the system (\ref{sys1}) satisfies Assumption \ref{ass2}. Then, the system (\ref{sys1}), (\ref{sys2}) is SP-UAS (UGAS) if the following conditions hold:
\begin{enumerate}
\item The system (\ref{sys1b}) is Lyapunov SP-UAS (Lyapunov UGAS);
\item The system (\ref{sys2}) is SP-UAS (UGAS);
\item The system (\ref{sys1}), (\ref{sys2}) satisfies the property USB (UGB). \\[-14mm]
\end{enumerate}
\null \hfill \null $\square$
\end{theorem}

For the sake of clarity we present the proofs of the main results in Section \ref{sec:proofs}. We provide only the proofs for the more general case of SP-UAS since the global versions follow by removing the restriction on the size of the domain of attraction and restricting the neighborhood of the origin to the origin itself (i.e., considering $\nu=0$).

\begin{rem}
While the proof of Theorem \ref{th1} is more direct since it does not require the existence of a Lyapunov function, Theorem \ref{th2} is 
{very} important since the existence of a Lyapunov function is 
{often helpful}  for controller design. Indeed, a Lyapunov function may allow us to improve the transients of the sampled-data system by redesigning a continuous-time controller in an appropriate way and using our results. 
\end{rem}

Two interesting corollaries for non-parameterized systems follow directly from our results. These results are interesting in cases when the exact discrete-time model of the plant can be computed and we do not have to appeal to Theorem \ref{th-motivate}. The appropriate definitions and assumptions for non-parameterized systems,
\begin{eqnarray}
x(k+1) & = & f(k,x(k),y(k)) \label{eq-nonpara-1} \\
y(k+1) & = & g(k,y(k)) \label{eq-nonpara-2}\,,
\end{eqnarray}
are obtained easily from global definitions for parameterized systems by setting $T=1$ and they are not repeated for space reasons. 

\begin{cor} Suppose that $f$ is continuous in $x,y$, uniformly in $k \geq 0$. Then, the system (\ref{eq-nonpara-1}), (\ref{eq-nonpara-2}) is UGAS if and only if:
\begin{enumerate}
\item The following system is UGAS
\begin{equation}
\label{eq-nonpara-3}
x(k+1)=f(k,x(k),0) \qquad ;
\end{equation}
\item The system (\ref{eq-nonpara-2}) is UGAS; and
\item The system (\ref{eq-nonpara-1}), (\ref{eq-nonpara-2}) is UGB.\\[-14mm]
\end{enumerate}
{    \null \hfill \null$\square$}
\end{cor} 

\begin{cor} Suppose that Assumption \ref{ass2} holds for $f$ in (\ref{eq-nonpara-1}). Then, the system (\ref{eq-nonpara-1}), (\ref{eq-nonpara-2}) is UGAS if:
\begin{enumerate}
\item The system (\ref{eq-nonpara-3}) is Lyapunov UGAS;
\item The system (\ref{eq-nonpara-2}) is UGAS; and
\item The system (\ref{eq-nonpara-1}), (\ref{eq-nonpara-2}) is UGB.\\[-14mm]
\end{enumerate}
{ \null \hfill \null$\square$}
\end{cor} 

We stress that UGB is in general difficult to check. {In \cite{DTCASCSCL} we present several sufficient conditions for this property to hold and which are inspired from \cite{MURAT-CASC,CASCAUT}. For the sake of completeness, we close this section with a result which is representative of the type of conditions given in those references and which we will use later on in the case study. To that end we introduce first the following technical hypothesis.
\begin{ass}\label{mars}
Suppose that there exist $\tilde{\alpha}_1,\tilde{\alpha}_2,\varphi \in {\cal K}_\infty$, {$\tilde{\gamma}_1$, $\tilde{\gamma}_2\in {\cal N}$}, $c,T^*>0$ and for each $T \in (0,T^*)$ there exists $V_T: \reals_{\geq 0} \times \reals^{n_x} \rightarrow \reals_{\geq 0}$ such that for all $x\in \reals^{n_x},z \in \reals^{n_z}$, $k \geq 0$ and $T \in (0,T^*)$ we have that :
\begin{eqnarray}
\tilde{\alpha}_1(\abs{x}) \ \leq \ V_T(k,x) & \leq & \tilde{\alpha}_2(\abs{x})+c \label{iisns-new-eq1} \\
V_T(k+1,f_T(k,x,z))-V_T(k+1,f_T(k,x,0)) & \leq & T \tilde{\gamma}_1(\abs{z}) \varphi(V_T(k,x))+ T\tilde{\gamma}_2(\abs{z}) \label{iisns-new-eq2} \\
V(k+1,f_T(k,x,0))-V_T(k,x) & \leq & 0  \label{iisns-eq3} \\
\int_1^\infty \frac{ds}{\varphi(s)} & = & \infty \ . \label{iisns-new-eq4} 
\end{eqnarray}
\end{ass}
\begin{prop}
\label{prop:iisns-lyap-new}
Consider the system (\ref{sys1}) with input $z$ and under Assumption \ref{mars}.  If furthermore the solutions of the system \rref{sys2} satisfy the summability condition
\begin{equation}
\label{eq:summability}
T \sum_{k=k_\circ}^{\infty} \mu \left( \abs{\phi_T^z(k,k_\circ,z_\circ)} \right) \leq \rho(\norm{x_\circ})\,,
\end{equation}
with some $\rho\in\cKinfty$ and $\mu(s):=\tilde\gamma_1(s) + \frac{\tilde\gamma_2(s)}{\varphi(1)}$ then, the system \rref{sys1} is UGB.
\null \hfill \null $\square$ 
\end{prop}
The proof follows along the lines of \cite[Lemmas 1 and 2]{DTCASCSCL} and therefore is omitted here\footnote{For review purposes we provide it in the Appendix.}. 

The proposition above establishes some interesting links with conditions used in the literature in the context of continuous-time systems to prove UGB. The condition (\ref{iisns-new-eq4}) restricts the growth of $\varphi$. In particular, it holds when $\varphi(s)=s$;  this situation was considered for instance  in \cite{CASCAUT} with  $\tilde{\gamma}_1(s) \equiv 0$. Earlier results using similar conditions are found in \cite{MAZPRA} and in the context of forward completeness already in \cite{SANCON}.

In the particular case when $\tilde \gamma_1$, $\tilde \gamma_2\in {\cal K}_\infty$ and $c=0$, one can prove that there exist $\cKinfty$ functions $\alpha_1$, $\alpha_2$  such that the systems trajectories satisfy 
\begin{eqnarray}\label{eq:iisns}
\alpha_1\left(\abs{\phi_T^x(k,k_\circ,x_\circ,\omega^z_{[k_\circ,k)})}\right) & \leq & \alpha_2(\abs{x_\circ})+ T \sum_{i=k_\circ}^{k-1} \mu(\abs{\omega^z_{[k_\circ,k)}}) , \qquad \forall k \geq k_\circ \geq 0 \,.
\end{eqnarray}
The property of the system by which it satisfies \rref{eq:iisns} is called integral  input-to-state neutral stability (iISNS) and has been introduced for {\em non}-parameterized discrete-time systems in \cite{ANG-dtIISS}.  See also \cite{IISS}, for the original definitions of integral input-to-state stability in the context of continuous-time systems. 

The summability condition \rref{eq:summability} which imposes a minimal convergence rate for the system \rref{sys2} is tight. However, we are not aware of any proof of necessity.  This condition is in general hard to check since it is trajectory-dependent. As far as we know such condition was firstly used in the context of continuous-time cascades, in \cite{CASCSCL}. Sufficient conditions for this property to hold are also given in \cite{DTCASCSCL} along the lines of the results presented in \cite{MURAT-CASC}.

Finally, we shall mention that the condition \rref{iisns-new-eq2} holds for instance (under suitable growth conditions on $V_T(k,\cdot)$) if $G_T(k,x,z)$ has linear growth in $x$ for each fixed $k$ and $z$.  This type of conditions have been extensively used in the continuous-time context to avoid the so called peaking phenomenon. See \cite{SEPBOOK} and references therein.
}

\section{A case study: tracking control of the unicycle}

\label{sec:appl}

In this section we revisit the problem of tracking control of a mobile robot of the unicycle type. This problem has been thoroughly studied in the continuous-time context via many different approaches (see \cite{KOLMCC} for a survey; for a more recent text with an updated list of references see \cite{ERJENSTHESIS}). To illustrate the utility of our results we will revisit the cascades approach used in \cite{MOBCAR} for a 3 degrees of freedom cart. The results may be extended to higher dimension systems, following for instance \cite{ERJENSTHESIS}. 

While the problem setting is the same as considered in the continuous-time context, we will see that the proof techniques employed in the discrete-time case are quite different. For instance, since we deal with approximate discrete-time models, some important structural characteristics are lost. Hence, we believe that the proofs of this section are interesting in their own right.

According to \cite{KANetal} the context of the problem can be set as follows. We have a mobile robot with two directional wheels and two ``fixed'' wheels and  whose motion is described by
\begin{equation}
\label{eq:system}
\dot x  =  v\cos\theta; \ \dot y  =  v\sin\theta; \ \dot\theta  =  \omega \ ,
\end{equation}
where $x$, $y$ are the Cartesian coordinates of the center of the axis joining the directional wheels and $\theta$ is the orientation angle of the directional wheels. The robot is required to follow a trajectory generated by an exosystem, i.e., a fictitious ``reference robot'' with kinematics 
\begin{equation}
\label{eq:system:ref}
\dot x_r  =  v_r(t)\cos\theta_r; \ \  \dot y_r  =  v_r(t)\sin\theta_r; \ \ \dot\theta_r  =  \omega_r(t)\,
\end{equation}
where $v_r(t)$ and $\omega_r(t)$ are $given$ reference velocities. Then, the tracking errors satisfy the set of equations (see \cite[Lemma 1]{KANetal})
\begin{equation}
\label{eq:errordynamics}
\dot x_e  =  \omega y_e-v+v_r(t)\cos\theta_e; \ \dot y_e  =  -\omega x_e+v_r(t)\sin\theta_e; \ \dot\theta_e  =  \omega_r(t)-\omega v \ ,
\end{equation}
where $(\cdot)_e := (\cdot )_r - (\cdot)$. The system is velocity-controlled, i.e., the control problem reduces to finding control inputs $\omega$ and $v$ (which also correspond to the actual angular and linear velocities of the cart) such that the origin of \rref{eq:errordynamics} is UGAS.

There are numerous solutions to this problem in the context of continuous-time (e.g. \cite{JIANONHOL96,SAMSIAM02} and \cite{ERJENSTHESIS} for a recent literature review). Here, we will revisit the cascaded approach proposed in \cite{MOBCAR} whose main feature is that the  control laws are $linear$. To illustrate and motivate our results we will proceed to solve the same problem with a $linear$ time-varying discrete-time controller which we will redesign based on the Euler-discretization of the error dynamics, 
\begin{equation}
\label{dt:eq:errordynamics}
\begin{array}{rcl}
x_e(k+1) & = & x_e(k) + T[\,\omega y_e(k)-v+v_r(k)\cos\theta_e(k)\,] \\
y_e(k+1) & = & y_e(k) + T[\,-\omega x_e(k)+v_r(k)\sin\theta_e(k)\,] \\
\theta_e(k+1) & = & \theta_e(k) + T[\,\omega_r(k)-\omega\,]\,.
\end{array}\,
\end{equation}
Thus, our control problem consists of designing $v$ and $\omega$ such that \rref{dt:eq:errordynamics} is UGAS.

We will solve this control problem following a similar approach to that of  \cite{MOBCAR} where it was shown using results for continuous-time cascaded systems, that the system \rref{eq:errordynamics} in closed loop with $v = v_r(t) + a_2x_e$ and $\omega= \omega_r(t) + a_1\theta_e$ is UGAS for appropriately chosen $a_1$ and $a_2$. 

The controller structure that we use
\begin{equation}\label{dt:control:mobcar}
\omega := \omega_r(k) + a_1\theta_e(k); \ \ v  := v_r(k) + a_2 x_e(k) + T\vartheta
\end{equation}
where $a_1$, $a_2$, $\omega_r(k)$ and $v_r(k)$ come from the continuous-time control law proposed in \cite{MOBCAR} and $\vartheta$ is an extra control input which depends on $k$, $x_e$ and $y_e$ and that we will design  with the aim of improving the system's performance. More specifically, the motivation for the control laws above is that as in the continuous-time context,  the closed loop system
\begin{eqnarray}\label{dt:cl:mobcar}
\label{dt:cl:mobcar:1}
&&
\begin{array}{ll}
 \begin{array}{rcl}
   x_e(k+1) & = &\\
   y_e(k+1) & = & 
 \end{array} 
\shiftleft{3mm}\underbrace{
 \begin{array}{l}
(1-Ta_2)x_e(k) + T\omega_r(k)y_e(k) - T^2\vartheta \\
y_e(k) - T\omega_r(k)x_e(k)
 \end{array}
 }_{\dty F_{1T}(k,x(k))} \ \, + & \\[12mm]
& \shiftleft{2in} \underbrace{
 \begin{array}{l}
T[\, a_1\theta_e(k)y_e(k) - v_r(k) + v_r(k)\cos\theta_e(k)  \,]\\
T[\, -a_1\theta_e(k)x_e(k) + v_r(k)\sin\theta_e(k)  \,]
 \end{array}
}_{\dty G_T(k,x(k),z(k))}
\end{array}\\
\label{dt:cl:mobcar:2}
&&\shiftright{2mm}
\begin{array}{rcl}
\theta_e(k+1) & = & (1-Ta_1)\theta_e(k) \ =: \  F_{2T}(k,x(k))\,,
\end{array} 
\end{eqnarray}
where $z:=\theta$ and $x:=\col[x_e,\, y_e]$, has a cascaded structure. 

Hence, the control laws \rref{dt:control:mobcar} are designed with two main ideas in mind: 1) to have as simple as possible controllers; 2) that the closed loop system verifies the conditions of our main results for cascades. More specifically, notice that the bottom subsystem \rref{dt:cl:mobcar:2} is independent of $x_e$ and $y_e$ and is UGES for  values of $a_1$ sufficiently small ($T^*a_1< 1$). Also, the interconnection term $G_T(k,x,z)$ is linear in $x_e$ and $y_e$. Hence, our results suggest that we only need to design $\vartheta$ as a function of $x_e$ and $y_e$ only so that the zero-input (i.e., with $G_T\equiv 0$) subsystem in \rref{dt:cl:mobcar:1} be UGAS (or possibly, UGES). 

Notice that in the particular case that  $\vartheta\equiv 0$ we obtain the emulated (discretized) continuous-time control law. However, as we will illustrate below, when carefully defined this extra degree of freedom in the control design allows to  improve performance and, on occasions, to enlarge the domain of attraction with respect to that of the emulated continuous-time control law. Thus, our control scheme can be regarded as a {\em redesign} of the cascaded-based continuous-time controller of \cite{MOBCAR}. Simulation results at the end of this section will illustrate the advantages of this approach. 

Here, we will establish UGES of the zero-input system based on a property of persistency of excitation. However, we will need a specific reformulation of this property\footnote{With respect to the original one from \cite{ASTBOH}.} within the framework of discrete-time parameterized systems. This is introduced next. To compact the notation,  in the sequel we will use $\omega_{r_k}:= \omega_r(k)$.
\begin{definition}[PE]
Let $\omega_r : \mZ_{\geq 0}\to\mZ $ be a function produced by sampling a function $\psi:\mR_{\geq 0}\to\mR$ at rate $T$. The function $\omega_r$ is said to be persistently exciting (PE) if there exist positive numbers $\mu$, $L$ and $T^*$ such that for all $T\in (0,T^*)$ and all $ j \geq  0$,
\begin{equation}
\label{pe:omegar}
T \sum_{k=j}^{j+\ell_{L,T}} \omega_{r_k}^2 \, \geq \, \mu\,. 
\end{equation}
\end{definition}

We are now ready to present our main result of this section.

\begin{prop}\label{prop:nonhol}
Consider the system \rref{dt:eq:errordynamics} in closed loop with \rref{dt:control:mobcar}. Assume that 
\begin{enumerate}
\item There exists $\hat{T},w_M > 0$ such that for all $k \geq 0$ and $T \in (0,\hat{T})$
\begin{equation} \label{wm}
\max\left\{ \abs{v_{r_k}}\,, \abs{\omega_{r_k}}, \frac{ \abs{\omega_{r_k}-\omega_{r_{k-1}}}}{T} \right\} \leq w_M \ .
\end{equation}
\item The signal $\omega_{r_k}$ is PE.
\end{enumerate}
Then, there exists $a_2>0$ such that for all $K,a_1>0$ and $\vartheta(k,x)$ with $\abs{\vartheta(k,x)} \leq K \abs{x}$, the system is UGAS.
\end{prop}
\noindent {\bf Proof. } It follows by invoking Theorem \ref{th2}. Firstly, we see that Assumption \ref{ass2} holds trivially in view of item 1 of the proposition. To see more clearly,  notice that \rref{wm} implies that $f_T(\cdot,\cdot,\cdot)$ is continuous and uniformly bounded in the first argument. Item 1 of the proposition also implies that there exists $c>0$ independent of $T$, such that for all $T\in (0,T^*)$, we have that $\norm{G_{T}(k,x,z)}\leq T c\norm{z}(\norm{x} + 1)$.

Secondly, it is evident that the origin of \rref{dt:cl:mobcar:2} is uniformly globally  exponentially stable for any $a_1$ and any $T \in (0,T^*)$ where $T^*>0$ is such that $1>a_1T^*>0$ and therefore the trajectories $\phi_T^z(\cdot,\cdot,\cdot)$ are uniformly summable. 

{It is left to prove that the unperturbed dynamics $x(k+1) = F_{1T}(k,x(k))$ in \rref{dt:cl:mobcar:1} is UGAS and that \rref{dt:cl:mobcar:1} is UGB. } As a matter of fact we will show that the zero-input system in \rref{dt:cl:mobcar:1} is Lyapunov UGES. 


\noindent {\underline{\sl Proof of UGES of $x(k+1) = F_{1T}(k,x(k))$}}: 
Consider the function $V_T(k,x):=\norm{x}^2 -\ep \omega_{r_{k-1}}x_ey_e$ with $\ep := \alpha_y + T$ and $\alpha_y >0$. Observe that this function  is positive definite and radially unbounded for sufficiently small $\alpha_y$, $T^*$ and $w_M$; indeed, we have that
\begin{equation}\label{c1c2V}
 c_1\norm{x}^2 \leq V_T(k,x) \leq c_2\norm{x}^2
\end{equation}
 with $c_1:=(1 - 0.5(\alpha_y+T^*)w_M)$ and $c_2:=(1+0.5(\alpha_y+T^*)w_M)$ which are clearly independent of $T$. We assume that $\alpha$ and $T^*$ are sufficiently small so that $c_1>0$. Next, we compute:
\begin{eqnarray}
\Delta V_T & := & V_T(k,F_{1T}(k,x)) - V_T(k,x)\nonumber\\ \label{deltaV}
& = &- T(2a_2 -\ep \omega_{r_k}^2) x_e^2 - \ep T \omega_{r_k}^2 y_e^2 \nonumber\\
&& + T^2 x_e \left( [a_2^2 + \omega_{r_k}^2 -\ep a_2 ] x_e - [2a_2 \omega_{r_k} - \ep \omega_{r_k}^3] y_e - [2(1-a_2 T) + \ep \omega_{r_k}^2 T ]\vartheta \right) \nonumber \\
&& + T^2 y_e \left(\omega_{r_k}^2y_e + \ep \omega_{r_k} \vartheta\right) + T^4\vartheta^2 + \ep x_ey_e (\omega_{r_k} - \omega_{r_{k-1}} + \omega_{r_k} a_2 T)\,.
\end{eqnarray}

Under the assumptions of the proposition, $\ep x_ey_e(\omega_{r_k} - \omega_{r_{k-1}} + \omega_{r_k} a_2 T) \leq (1/2)[T^2y_e^2 + \ep^2 w_M^2(1+a_2)^2x_e^2]$. Define $\alpha_x:=  a_2 -(\alpha_y+T^*)w_{M}^2 - (1/2)\ep^2 w_M^2(1+a_2)^2 $ which is positive for sufficiently small values of $\ep$ and sufficiently large values of $a_2$. Also, since $\abs{\vartheta(k,x)} \leq K\abs{x}$ there exists $K_1>0$ such that
\begin{equation}\label{93}
\frac{\Delta V_T}{T} \leq -(\alpha_x x_e^2 +  \alpha_y \omega_{r_k}^2 y_e^2) + T K_1 \norm{x}^2\, \ ,
\end{equation}
for all $T \in (0,T^*)$, $x \in \reals^2$ and $k\geq 0$. 

Let us introduce the following auxiliary function:
$$
W_T(k,x):= - T \sum^{\infty}_{i=k} e^{(k-i)T}\omega_{r_i} y_e^2 \ 
$$ 
for which we claim the following (for the proof, see the Appendix).
\begin{claim}
Suppose that
 the signal $\omega_{r_k}$ is PE. Suppose also that there exists $w_M>0$ such that for all $i \geq 0$ we have $\abs{\omega_{r_i}} \leq w_M$. Then, there exist strictly positive numbers $T^*,c_3,c_4,K_2,\tilde\alpha_y$ such that for all $T \in (0,T^*)$, $k \geq 0$ and $z_e\in\reals^2$ we have
\begin{eqnarray}
-c_3 y_e^2 & \leq & W_T(k,x) \ \leq \ -c_4 y_e^2 \label{w1} \\
\label{funfzig}\frac{\Delta W_T}{T} & \leq & \omega_{r_k}^2 y_e^2-\tilde\alpha_y y_e^2 + K_2 x_e^2 \label{w2} \ .
\end{eqnarray}
\ \hfill$\square$
\end{claim}

{Let now $T^*$ be generated by the claim above and such that (\ref{c1c2V}) and (\ref{93}) hold.  Then, we can } complete the proof 
by showing that there exists $\epsilon>0$ and $\tilde{T}>0$ such that $U_T(k,x):=V_T(k,x)+\epsilon W_T(k,x)$ is a Lyapunov function that proves UGES. Indeed, let   $\epsilon: = \min \left\{ \frac{c_1}{2{c_3}}, \frac{\alpha_x}{2 K_2}, {\alpha}_y \right\}$ and $\tilde{T} := \min \left\{ T^* , {\frac{1}{2K_1}} \min\left\{ \frac{\alpha_x}{2}, \epsilon {\tilde\alpha_y} \right\} \right\}$, where $T^*>0$ comes from the proof of the proposition {and $\tilde\alpha_y$ comes from \rref{funfzig}}. Then it is easy to show that $U_T$ satisfies
\begin{eqnarray}
\label{eq:U1}
\frac{c_1}{2} \abs{x}^2& \leq & U_T(k,x) \ \leq c_2 \abs{x}^2 \\
\label{eq:U2} 
\frac{\Delta U_T}{T} & \leq & -\tilde{c}_3 \abs{x}^2 \ , 
\end{eqnarray}
where $\tilde{c}_3=\frac{1}{2}\min\{ \frac{\alpha_x}{2},\epsilon {\tilde\alpha_y} \}$. This completes the proof {invoking standard Lyapunov arguments}.

\noindent {\underline{\sl Proof of UGB}: We invoke Proposition \ref{prop:iisns-lyap-new}. It is worth recalling to avoid confusion in the notation, that the state $x$ in Proposition   \ref{prop:iisns-lyap-new} corresponds here to $x = \col[x_e,\,y_e]$ and the input $z$ in Proposition \ref{prop:iisns-lyap-new} corresponds here to $\theta_e$. Hence, we proceed to verify the conditions of the proposition with $f_T(k,x,z):=F_{1T}(k,x)+G_T(k,x,z)$ as defined in \rref{dt:cl:mobcar:1} and  $V_T(k,x) = U_T(k,x)$. The bounds \rref{iisns-new-eq1} and \rref{iisns-eq3} hold from \rref{eq:U1} and \rref{eq:U2}. The conditions \rref{iisns-new-eq2}  and  \rref{iisns-new-eq4} hold with $\varphi(s)=s$, $\tilde\gamma_2 = d_2 s$, and $\tilde\gamma_1(s) := d_3 s$, $d_2,\, d_3>0$.  This is because $U_T(k,\cdot)$ is quadratic and $G_T(k,x,z)$ contains terms of linear growth in $x$  for each fixed $k$ and $z$ and trigonometric functions, which can be over-bounded by a linear function of $x$. Also, $G_T(k,x,\cdot)$ can be over-bounded by a linear function for each fixed $k$ and $x$. Finally, \rref{eq:summability} holds with a linear function $\rho(s):=\rho\, s$ since $\mu(s)$ is in this case a non decreasing function of linear growth and $\phi_T^{\theta_e}(k)$ decays uniformly exponentially to zero. 
\null\hfill\null $\blacksquare$}

\begin{rem}\label{rem:proofUGC:cs}
It is worth remarking that in the proof above we could also have invoked Theorem \ref{th1}. Notice that to show that the system \rref{sys1} is UGC with input $\theta_e$ we may appeal to Lemma \ref{lem3} observing that the system $x(k+1) = f_T(k,x(k),\theta_e(k))$ with $f_T(k,x,\theta_e):= F_{1T}(k,x) + G_{1T}(k,x,\theta_e)$ is linear in the state $x=\col[x_e,\,y_e]$ and for each $\Delta>0$, we have that $f_T(k,x,\theta_e)$ is globally Lipschitz in $\theta_e$ uniformly for all $k \geq 0$ and $x$ such that $\norm{x}\leq \Delta$, with a Lipschitz constant of the form $L= L'T$ with $L'$ depending only on $\Delta$. 
\end{rem}

\putfig{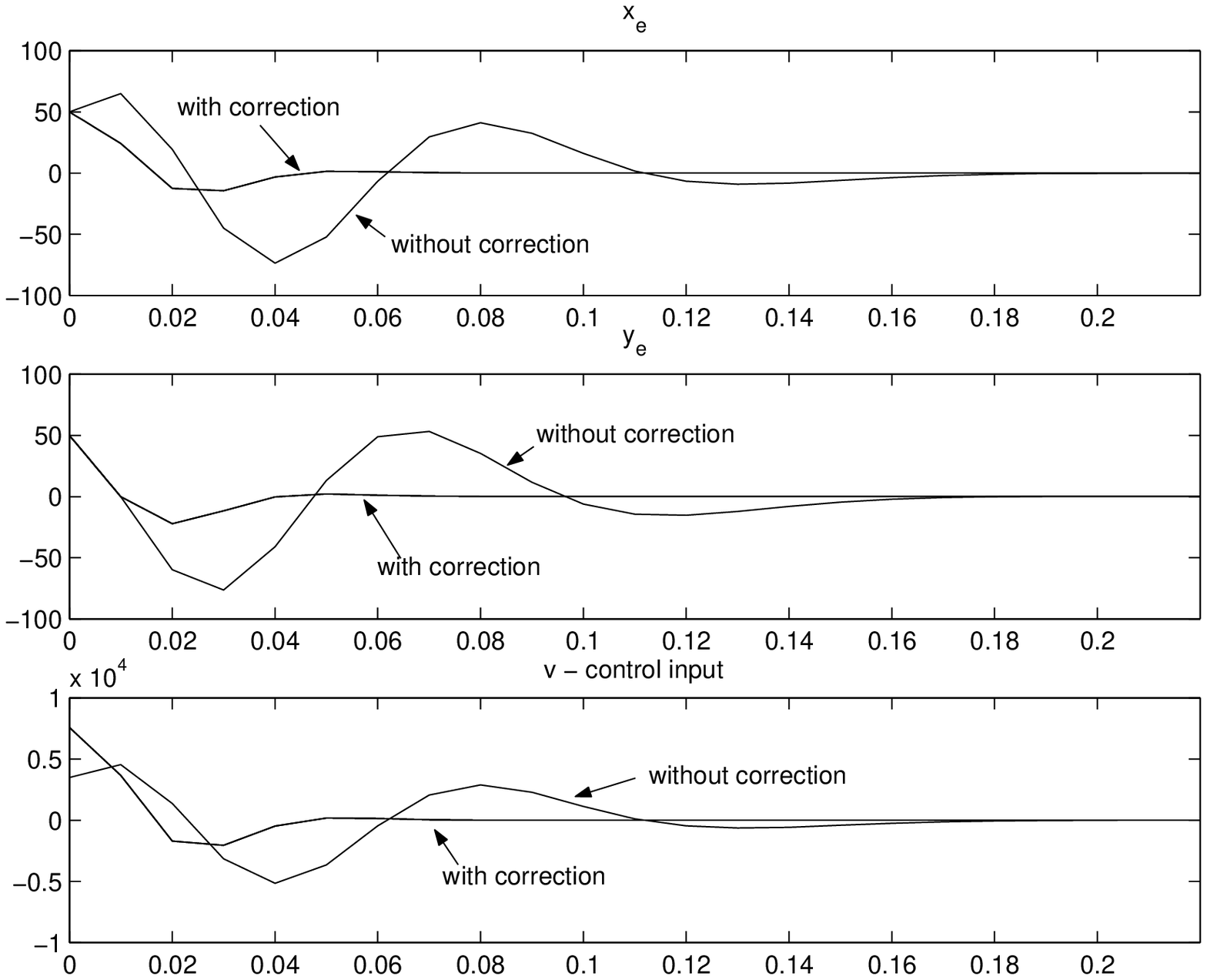}{simu3_nice}{70}{Tracking errors for $x$ and $y$ and the control input $v$ given by \rref{dt:control:mobcar}, \rref{correction:nonhol}.}

Now we illustrate how we can use Proposition \ref{prop:nonhol} to improve the performance of the system with the redesigned controller. Since the correction $\vartheta$ can be chosen arbitrarily, we can choose it so that negativity of $\Delta V$ in \rref{deltaV} is enhanced. Indeed, a closer inspection of the difference equation \rref{deltaV} shows that one such choice is 
\begin{equation}\label{correction:nonhol}
\vartheta(k,x) := \dty\frac{(a_2^2 + \omega_{r_k}^2 - \ep a_2) x_e - (2a_2 \omega_{r_k} - \ep \omega_{r_k}^3 ) y_e}{2(1-a_2T)+\ep\omega_{r_k}^2T}
\end{equation}
with $\ep = \alpha_y + T$. 

The simulations that we present next show that the performance is considerably improved. We have simulated the system above in \simulink\ of \matlab\ with $a_2= 70$, $a_1=10$, $w_r(k)=20\sin(kT)$, $T = 0.01$ and $\alpha_y = 2-T$. We show the results in Figure \ref{simu3_nice}. We show only the responses for the states $x_e$ and $y_e$ as well as $v$ since these are the only variables affected by the additional input $T\vartheta$. We show simulations for the system's response with $\vartheta = 0$, with $\vartheta = 0.5[\,(a_2^2 + \omega_{r_k}^2 - \ep a_2) x_e - (2a_2 \omega_{r_k} - \ep \omega_{r_k}^3 ) y_e\,]$ and for $\vartheta$ as defined in \rref{correction:nonhol}.  The best apparent performance is for the latter. 

It is also clear from the plots, that even though the correction $\vartheta$ is linear in the state  and actually ($\vartheta \approx {\cal O}(1)\norm{x}$), this correction is not comparable to ``adding gain'' to the control input. 
Notice that in this case the resulting control effort is actually smaller than in the case of the continuous-time based controller (i.e., when $\vartheta=0$). 

\section{Proofs of the main theorems}
\label{sec:proofs}
All results are only proved for semiglobal-practical stability properties since global results follow the same steps with minor changes.

\subsection{Proof of Theorem \ref{th1}}
\begin{lemma} 
\lab{lem1}
Suppose that the solutions of the system (\ref{sys1}) with input $z$ are USC (UGC) and the system (\ref{sys1b}) is SP-UAS (UGAS), with the function $\beta_x \in {\cal KL}$. Then, for any pair of strictly positive real numbers $(\Delta,\nu)$ there exists $T^*>0$ (there exists $T^*$) such that for any $\eta \in (0,\Delta]$ (any $\eta >0$) and any strictly positive real numbers $(\epsilon,L)$ there exists $\mu>0$ such that for all inputs $z(\cdot)$ with $\norm{\omega^z} \leq \mu$, $k_\circ \geq 0$, $x(k_\circ)=x_\circ$ with $\abs{x_\circ} \leq \eta$ and $T \in (0,T^*)$ the following holds:
\begin{equation}
\label{eq:klstab}
\abs{\phi^x_T(k,k_\circ,x_\circ,\omega^z_{[k_\circ,k)})} \leq \max\{ \beta_x(\abs{x_\circ},(k-k_\circ)T), \nu \}+\epsilon, 
\end{equation}
for all $k \in [k_\circ,k_\circ+\ell_{L,T}]$.
\null \hfill \null $\square$
\end{lemma}

\noindent {\bf Proof. } 
For any $(\Delta,\nu)$ there exists $T^*>0$ such that for any $\eta \in (0,\Delta]$ and any strictly positive real numbers $(\epsilon,L)$ there exists $\mu>0$ such that for all inputs $z(\cdot)$ with $\norm{\omega^z} \leq \mu$, all $k_\circ\geq 0$, $x(k_\circ)=x_\circ$ with $\abs{x_\circ} \leq \eta$ and $T \in (0,T^*)$ the following holds:
\begin{eqnarray}
\abs{\phi^x_T(k,k_\circ,x_\circ,\omega^z_{[k_\circ,k)})} & \leq & \abs{\phi^x_T(k,k_\circ,x_\circ,0)}+\abs{\phi^x_T(k,k_\circ,x_\circ,\omega^z_{[k_\circ,k)})-\phi^x_T(k,k_\circ,x_\circ,0)} \nonumber 
\end{eqnarray}
Hence, \rref{eq:klstab} follows from \rref{gato} and the ${\cal K}{\cal L}$ bound on $\abs{\phi^x_T(k,k_\circ,x_\circ,0)} $.

\begin{lemma}
\label{cor1}
Suppose that all the conditions of Lemma \ref{lem1} hold. Then, for any strictly positive $(\Delta,\nu)$ there exists $T^*>0$ such that for any $\epsilon>0$ and $\eta \in (0,\Delta]$ there exists $\mu>0$ such that for any input $z(\cdot)$ with $\norm{\omega^z} < \mu$, $k_\circ \geq 0$, $x(k_\circ)=x_\circ$ with $\abs{x_\circ} \leq \Delta$ and $T \in (0,T^*)$ we have that the solutions of the system (\ref{sys1}) with input $z$ satisfy the inequality (\ref{eq:klstab}) for all $k \geq k_\circ$.
\null \hfill \null $\square$
\end{lemma}

\noindent {\bf Proof. } 
Let $\beta_x \in {\cal KL}$ come from SP-UAS stability of (\ref{sys1b}). Let $(\Delta,\nu)$ be given. Without loss of generality assume that $\Delta>\nu$ and 
\begin{equation}
\label{beta0}
\beta_x(s,0) \geq s, \qquad\forall s \geq 0\,.
\end{equation}
 Let $\tilde{\nu} \in (0,1)$ be such that 
\begin{equation}
\label{tildenu}
\beta_x(\tilde{\nu},0)+\tilde{\nu}/2 \leq \nu \ .
\end{equation}
Note that (\ref{beta0}) implies that $\tilde{\nu}<\nu$. Let $(\Delta,\tilde{\nu}/2)$ generate $T^*>0$ via Lemma \ref{lem1} and assume without loss of generality that $T^*<1$. Let $\epsilon>0$ and $\eta \in (0,\Delta]$ be given. Let $\delta_1>0$ be such that 
\begin{equation}
\label{delta1}
\beta_x(\delta_1,0)+\delta_1/2 \leq \epsilon \ .
\end{equation}
Let $L>1$ be such that 
\begin{equation}
\label{L}
\beta_x(\max\{ \eta,\delta_1, 1\},L-1)<\frac{\delta_1}{2} \ .
\end{equation}
Let $\epsilon_1:=\min\{ \tilde{\nu}/2,\delta_1/2\}$ and let $\epsilon_1$ and $L$ generate $\mu>0$ via Lemma \ref{lem1}. In the rest of the proof we consider arbitrary fixed $k_\circ \geq 0$, $x(k_\circ)=x_\circ$ with $\abs{x_\circ} \leq \eta$, $\norm{\omega^z} \leq \mu$ and $T \in (0,T^*)$. 

In order to simplify the notation, for $i\in\mZ_{\geq 1}$ we introduce:
\begin{eqnarray}
k_i & := & k_\circ+ i \ell_{L,T} \nonumber \\
x_i & := & \phi_T^x(k_i,k_\circ,x_\circ,\omega^z_{[k_\circ,k_i)}) \nonumber \ .
\end{eqnarray}

\noindent{\bf Time interval $[k_\circ,k_1]$:} From Lemma \ref{lem1}, (\ref{tildenu}), (\ref{beta0}) 
we can write
\begin{eqnarray}
\label{step:1a}
\abs{\phi^x_T(k,k_\circ,x_\circ,\omega^z_{[k_\circ,k)})} & \leq & \max\{ \beta_x(\abs{x_\circ},(k-k_\circ)T), \tilde{\nu}/2 \}+\epsilon_1 \nonumber \\
  & \leq & \max\{ \beta_x(\abs{x_\circ},(k-k_\circ)T), \nu \}+ \epsilon
\end{eqnarray}
for all $k \in [k_\circ,k_1]$. Moreover, note that
\begin{equation}
\label{L1}
T<1 \qquad \Longrightarrow \qquad L-1 \leq T\ell_{L,T} \ . 
\end{equation}
Using, (\ref{L}), (\ref{L1}) and our choice of $\epsilon_1$ we can write
\begin{eqnarray}
\label{step:1b}
\abs{x_1} \ = \ \abs{\phi^x_T(k_1,k_\circ,x_\circ,\omega^z_{[k_\circ,k_1)})} & \leq & \max\{ \beta_x(\abs{x_\circ}, T \ell_{L,T}), \tilde{\nu}/2 \}+\epsilon_1 \nonumber \\
 & \leq & \max\{ \beta_x(\eta,L-1), \tilde{\nu}/2 \}+\epsilon_1 \nonumber \\
 & \leq & \max\{\delta_1/2,\tilde{\nu}/2\}+\min\{\delta_1/2,\tilde{\nu}/2 \} \nonumber \\
 & \leq & \max\{\delta_1,\tilde{\nu} \} \ .
\end{eqnarray}
 
\noindent{\bf Time interval $[k_1,k_2]$:}  Since $\abs{x_1} \leq \max\{ \delta_1,\tilde{\nu} \}$, we consider two cases. \\
{\em Case 1:} $\abs{x_1} \leq \tilde{\nu}$. In this case, using Lemma \ref{lem1}, (\ref{tildenu}) and (\ref{beta0}) we can write
\begin{eqnarray}
\label{step2:c1a}
\abs{\phi_T^x(k,k_\circ,x_\circ,\omega^z_{[k_\circ,k)})} & = & \abs{\phi^x_T(k,k_1,x_1,\omega^z_{[k_1,k)})} \nonumber \\
 & \leq & \max \{ \beta_x(\abs{x_1}, 0),\tilde{\nu}/2 \}+\epsilon_1 \nonumber \\
  & \leq & \max\{ \beta_x(\tilde{\nu},0), \tilde{\nu}/2 \} +\min\{ \tilde{\nu}/2,\delta_1/2 \}  \\
  & \leq & \beta_x(\tilde{\nu},0)+\tilde{\nu}/2 \ \leq \ \nu \ , \nonumber 
\end{eqnarray}
for all $k \in [k_1,k_2]$. Moreover, using (\ref{step:1b}), the fact that $\tilde{\nu}<1$, (\ref{L}), (\ref{L1}) and our choice of $\epsilon_1$ we can write
\begin{eqnarray}
\label{step2:c1b}
\abs{\phi^x_T(k_2,k_\circ,x_\circ,\omega^z_{[k_\circ,k_2)})} & \leq & \max \{ \beta_x(\abs{x_1}, T\ell_{L,T}),\tilde{\nu}/2 \}+\epsilon_1 \nonumber \\
 & \leq & \max \{ \beta_x(\tilde{\nu}, L-1),\tilde{\nu}/2 \}+\epsilon_1 \nonumber \\
 & \leq & \max \{ \beta_x(1, L-1),\tilde{\nu}/2 \}+\epsilon_1  \\
 & \leq & \max\{ \delta_1/2,\tilde{\nu}/2 \}+\min\{\delta_1/2,\tilde{\nu}/2 \} \ \leq \ \max\{ \delta_1,\tilde{\nu} \} \ . \nonumber
\end{eqnarray}
{\em Case 2:} $\abs{x_1} \leq \delta_1$. Using Lemma \ref{lem1} and (\ref{delta1}) we can write
\begin{eqnarray}
\label{step2:c2a}
\abs{\phi^x_T(k,k_\circ,x_\circ,\omega^z_{[k_\circ,k)})} & = & \abs{\phi^x_T(k,k_1,x_1,\omega^z_{[k_1,k)})} \nonumber \\
  & \leq & \max\{ \beta_x(\delta_1,0), \tilde{\nu}/2 \} +\min\{ \tilde{\nu}/2,\delta_1/2 \}  \\
  & \leq & \max\{ \beta_x(\delta_1,0)+\delta_1/2,\tilde{\nu} \} \ \leq \ \max\{ \epsilon, \tilde{\nu} \} \ ,
\end{eqnarray}
for all $k \in [k_1,k_2]$. Also, in a similar manner as before it follows that
\begin{eqnarray}
\label{step2:c2b}
\abs{\phi^x_T(k_2,k_\circ,x_\circ,\omega^z_{[k_\circ,k_2)})} & = & \abs{\phi^x_T(k_2,k_1,x_1,\omega^z_{[k_1,k_2)})} \nonumber \\
  & \leq & \max\{ \beta_x(\delta_1,T \ell_{L,T}), \tilde{\nu}/2 \} + \epsilon_1  \\
  & \leq & \max\{ \beta_x(\delta_1,L-1), \tilde{\nu}/2 \} +\epsilon_1  \\
  & \leq & \max\{ \delta_1/2,\tilde{\nu} \} +\min\{ \tilde{\nu}/2,\delta_1/2 \} \nonumber \\
 & \leq & \max\{ \delta_1, \tilde{\nu} \} \ .
\end{eqnarray}

\noindent{\bf Time intervals $[k_i,k_{i+1}]$, $i\geq 1$:} Using similar calculations it can be shown by induction that for all integers $i \geq 1$,
\begin{eqnarray}
\abs{\phi^x_T(k,k_\circ,x_\circ,\omega^z_{[k_\circ,k)})} & \leq & \max\{ \epsilon, \nu \} \ , \ \forall k \in [k_i,k_{i+1}] \label{allsteps} \\
\abs{\phi^x_T(k_i,k_\circ,x_\circ,\omega^z_{[k_\circ,k_i)})}  & \leq  & \max\{ \delta_1, \tilde{\nu} \} \ .
\end{eqnarray}
The proof follows from (\ref{step:1a}) and (\ref{allsteps}) by noting that
\begin{eqnarray}
\abs{\phi^x_T(k,k_\circ,x_\circ,\omega^z_{[k_\circ,k)})} & \leq & \max\{ \epsilon, \nu \} \ \leq \ \epsilon+\nu \nonumber \\
 & \leq & \max\{ \beta_x(\abs{x_\circ},(k-k_\circ)T),\nu\}+\epsilon \ , \forall k \geq k_1 \ .
\end{eqnarray} 

\begin{lemma}
\label{lem2}
Suppose that all conditions of Lemma \ref{lem1} hold and $\beta_x \in {\cal KL}$ comes from SP-UAS (UGAS) of (\ref{sys1b}). Then, there exist $c_x,c_u>0$, $\gamma \in {\cal K}$ and for any $\nu>0$ there exists $T^*>0$ such that for all $k_\circ \geq 0$, $x(k_\circ)=x_\circ$ with $\abs{x_\circ} \leq c_x$, $\norm{\omega^z} < c_u$ and $T \in (0,T^*)$ the following holds:
\begin{equation}\label{33}
\abs{\phi^x_T(k,k_\circ,x_\circ,\omega^z_{[k_\circ,k)})}  \leq \max\{ \beta_x(\abs{x_\circ},(k-k_\circ)T) , \nu \} + \gamma\left(\norm{\omega^z_{[k_\circ,k)}}\right) \ 
\end{equation}
for all $k \geq k_\circ \geq 0$. \null \hfill \null $\square$
\end{lemma}

\noindent {\bf Proof. } 
The proof of this fact follows closely the proof of \cite[Lemma 3.3]{KHALIL}. Since all conditions of Lemma \ref{lem1} hold, the conclusion of Lemma \ref{cor1} holds. Let $(c_x,\nu)$ generate $T^*>0$ via Lemma \ref{lem1}. Let $\eta=c_x$ and for a fixed $\epsilon>0$, let $\bar{\mu}(\epsilon)$ be the supremum over all applicable $\mu$'s. Then, $\norm{\omega^z_{[k_\circ,k)}} < \bar{\mu}(\epsilon)$ implies that (\ref{eq:klstab}) holds for all $k \geq k_\circ \geq 0$ and if $\mu_1 >\bar{\mu}(\epsilon)$ then there exists $\norm{\bar{\omega}^z_{[k_\circ,k)}}<\mu_1$ and $\abs{\phi^x_T(k,k_\circ,x_\circ,\bar{\omega}^z_{[k_\circ,k)})} > \max\{ \beta_x(\abs{x_\circ},(k-k_\circ)T),\nu\}+\epsilon$. $\bar{\mu}(\epsilon)$ is positive and nondecreasing but it is not necessarily continuous. Chose $\zeta \in {\cal K}$ such that $\zeta(r) \leq k\bar{\mu}(r)$, with $k \in (0,1)$. Let $\gamma(r)=\zeta^{-1}(r)$ and note that $\gamma \in {\cal K}$. Let $c_u:=\lim_{r \rightarrow \infty} \zeta(r)$. 

Let $\abs{x_\circ} \leq c_x$, $T \in (0,T^*)$. Given $\omega^z_{[k_\circ,k)}$ with $\norm{\omega^z_{[k_\circ,k)}}<c_u$ let $\epsilon=\gamma\left(\norm{\omega^z_{[k_\circ,k)}}\right)$. Then, we have that \rref{33} holds 
for all $k \geq k_\circ \geq 0$.

\begin{rem} \lab{betas}
The following fact was proved in \cite{NESTEESON}: for any $\beta \in {\cal KL}$ and any $c>0$ there exists $\tilde{\beta} \in {\cal KL}$ such that
$$
\beta(s,t) \leq \tilde{\beta}(s,t+c)\,, \quad \forall s,t \in [0,\infty) \ .
$$
This can be further strengthened in the following manner. It was shown in \cite{IISS} that given any $\beta \in {\cal KL}$, there exist $\sigma\in {\cal K}_\infty$, $\kappa\in {\cal K}$ such that $\beta(r,t) \leq \sigma(\kappa(r)e^{-t})$ $\forall\, r,t \geq 0$. Consequently, given any $\beta \in {\cal KL}$ and a nondecreasing function $L(\cdot)$, there exists $\tilde{\beta} \in {\cal KL}$ such that 
$$
\beta(s,t) \leq \tilde{\beta}(s,t), \qquad \forall t \in [0,L(s)] \ .
$$
To show this, let $\beta$ generate $\sigma,\kappa \in {\cal K}_\infty$ as above. Then $\tilde{\beta}(s,t):= \sigma(\kappa(s)e^{L(s)}e^{-t})$ proves the claim.
\null \hfill \null $\square$
\end{rem}

\begin{lemma} 
\label{lem-kl}
Let $\xi=(x^T \ z^T)^T$. Suppose that there exist $\beta_1,\beta_2,\beta_3 \in {\cal KL}$, $\gamma \in {\cal K}$ and $c_0>0$ such that for all $\nu_1>0$, there exists $T^*>0$ such that for all $\max\{\abs{u(k_\circ)}, \abs{\xi(k_\circ)}\} \leq c_0$, $T \in (0,T^*)$ and $k \geq k_\circ \geq 0$ we have

\begin{eqnarray}
\abs{x(k)} & \leq & \max\{ \beta_1(\abs{x(k_\circ)},(k-k_\circ)T),\nu_1\}+\gamma \left(\sup_{t \in [k_\circ,k]} \abs{u(t)} \right) \nonumber \\
\abs{u(k)} & \leq & \max\{ \beta_2(\abs{\xi(k_\circ)},(k-k_\circ)T),\nu_1\} \ , \nonumber \\
\abs{z(k)} & \leq & \max\{ \beta_3(\abs{z(k_\circ)},(k-k_\circ)T),\nu_1\} \ ,
\end{eqnarray}
where $\xi:=(x^T \ z^T)^T$. Then, there exist $\beta \in {\cal KL}$ and $c_1>0$ such that for any $\nu>0$ there exists $T^*_1>0$ such that for all $\xi(k_\circ)=\xi_\circ$ with $\abs{\xi(k_\circ)} \leq c_1$, $T\in (0,T^*)$ and $k \geq k_\circ \geq 0$ we have
$$
\abs{\xi(k)} \leq \max\{ \beta(\abs{\xi_\circ}, (k-k_\circ)T),\nu \} \ . 
$$ 
In particular, we can take:
\begin{eqnarray}
\beta(s,t) & = & 4\tilde{\beta}_1\left( 2 \tilde{\beta}_1(s,t/2)+2\gamma(\tilde{\beta}_2(s,0)), t/2 \right)+4 \gamma(\tilde{\beta}_2(s,t/2))+2 \tilde{\beta}_3(s,t) \ , \label{new-beta}
\end{eqnarray}
where all $\tilde{\beta}_i$ are generated via $\beta_i$ and an arbitrary fixed $c>0$ using Remark \ref{betas}. Moreover, if all the inequalities hold globally, then the conclusion holds globally and we can take:
\begin{eqnarray}
\beta(s,t) & = & \tilde{\beta}_1\left( \tilde{\beta}_1(s,t/2)+\gamma(\tilde{\beta}_2(s,0)), t/2 \right)+ \gamma(\tilde{\beta}_2(s,t/2))+\tilde{\beta}_3(s,t) \ .
\end{eqnarray}
\null \hfill \null $\square$
\end{lemma}

\noindent {\bf Proof. } 
Let $c>0$ be fixed and let $\beta_i$ and $c$ generate $\tilde{\beta}_i$ via Remark \ref{betas}. Let $c_0$ come from the Lemma and let $\nu$ be given.  Let $c_1,\nu_1$ be strictly positive numbers such that
\begin{eqnarray}
c_0 & \geq & \max\{ 2\tilde{\beta}_1(c_1,0)+2\gamma(\tilde{\beta}_1(c_1,0)), \tilde{\beta}_2(c_1,0),c_1 \} \nonumber \\
\nu & \geq & \gamma_0(\nu_1) \ ,
\end{eqnarray}
where $\gamma_0(s):=4 \tilde{\beta}_1(2s+2\gamma(s),0)+4 \gamma(s)+2s$.
Let $\nu_1$ generate $T_1^*>0$ via the Lemma. Let $T^*=\min\{ T^*_1, c \}$. Define a piecewise constant function $x(\cdot)$ as $x(t)=x(k), t \in [kT,(k+1)T)$;  $z(\cdot),u(\cdot),\xi(\cdot)$ are defined in the same manner and can be regarded as functions that depend on continuous time $t \in \reals$. From the definition of $\tilde{\beta}_i$ and our choice of $T^*$ we have $\tilde{\beta}_i(s,t+T^*) \geq \tilde{\beta}_i(s,t+c) \geq  \beta_i(s,t)$, $\forall s,t$. 

Let $\abs{\xi(t_\circ)} \leq c_1$ and $T \in (0,T^*)$. Then the following holds:
\begin{eqnarray}
\abs{x(t)} & \leq & \max\{ \tilde{\beta}_1(\abs{x(t_\circ)},t-t_\circ),\nu_1\}+\gamma \left(\sup_{\tau \in [t_\circ,t]} \abs{u(\tau)} \right) \nonumber \\
\abs{u(t)} & \leq & \max\{ \tilde{\beta}_2(\abs{\xi(t_\circ)},t-t_\circ),\nu_1\} \ , \nonumber \\
\abs{z(t)} & \leq & \max\{ \tilde{\beta}_3(\abs{z(t_\circ)},t-t_\circ),\nu_1\} \ ,
\end{eqnarray}
for all $t \geq t_\circ \geq 0$. Using the same technique as in \cite[pg. 221-222]{KHALIL}, we first find a bound for $x(t)$ starting at $\frac{t+t_\circ}{2}$ and then a bound for $x\left(\frac{t+t_\circ}{2}\right)$ starting at $t_\circ$. Using the same procedure as in \cite[pg. 221-222]{KHALIL}, these inequalities together with the following facts
\begin{eqnarray}
\max\{ a,b \}+\max\{ c,d \} & \leq & \max\{ 2a+2c,2b+2d\} \ ,\quad \forall a,b,c,d \geq 0 \nonumber \\
\sup_{\tau \in [t_\circ,(t+t_\circ)/2]} \abs{u(\tau)} & \leq & \tilde{\beta}_2(\abs{\xi(t_\circ)},0) \nonumber \\
\sup_{\tau \in [(t+t_\circ)/2,t]} \abs{u(\tau)} & \leq & \tilde{\beta}_2(\abs{\xi(t_\circ)},(t-t_\circ)/2) \nonumber \\
\abs{\xi(t)} & \leq & \abs{x(t)} + \abs{z(t)}\,, \nonumber
\end{eqnarray}
generate the inequality
$$
\abs{\xi(t)} \leq \max\{ \beta(\abs{\xi(t_\circ)},t-t_\circ),\gamma_0(\nu_1) \}\,.
$$
This inequality, together with our choice of $\nu_1$ and the fact that $\xi(kT)=\xi(k)$ implies that
$$
\abs{\xi(k)} \leq \max\{ \beta(\abs{\xi(k_\circ)},(k-k_\circ)T), \nu \} \ 
$$
which completes the proof.

\noindent {\bf Proof of Theorem \ref{th1}:} 
Let $\beta_x \in {\cal KL}$ come from item 1 of the Theorem. Let $\beta_z \in {\cal KL}$ come from item 2 of the Theorem. Let $\kappa \in {\cal K}_\infty$ and $c>0$ come from item 3 of the Theorem. Let $\gamma \in {\cal K}$ and $c_x,c_u$ come from Lemma \ref{lem2}. Let $c_z>0$ be such that $\beta_z(c_z,0) \leq c_u$ and $c_\xi:=\min\{ c_z,c_x \}$. Let $\bar{\beta} \in {\cal KL}$ and $c_1>0$ be generated using $\beta_1=\beta_x,\beta_2=\beta_3=\beta_z,\gamma$ and $c_0=c_\xi$ via Lemma \ref{lem-kl}. Let $\kappa_1(s):=\left(\frac{c}{\kappa(c_1)}+1\right) \kappa(s)$. Let $L_i:\reals_{\geq 0} \rightarrow \reals_{\geq 0}, i=1,2$ be continuous nondecreasing functions with $L_i(0)= 1$ such that for all $\Delta>0$ we have
\begin{eqnarray}
\beta_z(\Delta,L_1(\Delta)-1) & \leq & c_1  \lab{betaz} \\
\beta_x(\kappa_1(\Delta),L_2 \circ \kappa_1(\Delta)-1) & \leq & c_1/2 \lab{betax} \,.
\end{eqnarray} 
Let $L(s):=L_1(s)+L_2 \circ \kappa_1(s)$. Let $\bar{\beta}$ and $L(s)$ generate $\tilde{\beta}$ via Remark \ref{betas}.
Finally, we define
$$
\beta(s,t):=\max\{ \tilde{\beta}(\kappa_1(s),t),\bar{\beta}(s,t),\kappa_1(s)e^{L_1(s)+L_2 \circ \kappa_1(s)}e^{-t} \} \ .
$$
Let $(\Delta,\nu)$ be given. Let $\nu_1:=\min\{ c_1/2,\nu\}$. Let $(\kappa_1(\Delta),\nu_1)$ generate $T_1^*>0$ via Lemma \ref{lem1}; let $(\Delta,\nu_1)$ generate $T_2^*>0$ via the item 2 of the Theorem; let $\nu_1$ generate $T_3^*>0$ via Lemma \ref{lem2}; let $\nu_1>0$ generate $T_4^*>0$ via Lemma \ref{lem-kl}. Let $T^*:=\min\{ T_1^*,T_2^*,T_3^*,T_4^*,1\}$ and $T \in (0,T^*)$.  Let $k_\circ \geq 0$, $\xi(k_\circ)=\xi_\circ$ with $\abs{\xi_\circ} \leq \Delta$. We consider two cases.

{\em Case 1:} If $\abs{\xi_\circ} \leq c_1$, then by a direct application of Lemma \ref{lem2} and \ref{lem-kl} we have that
\begin{equation}
\label{case1-bound}
\abs{\phi_T^\xi(k,k_\circ,\xi_\circ)} \leq \max\{ \bar{\beta}(\abs{\xi_\circ},(k-k_\circ)T),\nu_1 \} \leq \max\{ \bar{\beta}(\abs{\xi_\circ},(k-k_\circ)T),\nu \} ,
\end{equation}
for all $k \geq k_\circ \geq 0$.

{\em Case 2:} If $\abs{\xi_\circ} \in [c_1,\Delta]$, then since $\abs{\xi_\circ} \geq c_1$, we can write using item 3 of the Theorem:
\begin{eqnarray}
\abs{\phi_T^\xi(k,k_\circ,\xi_\circ)} & \leq & \kappa(\abs{\xi_\circ})+c \ \leq \ \kappa(\abs{\xi_\circ})+\frac{c}{\kappa(c_1)}\kappa(\abs{\xi_\circ})=\kappa_1(\abs{\xi_\circ}) \ ,
\end{eqnarray}
for all $k \geq k_\circ \geq 0$. Hence, we can write that for all $k \in [k_\circ,k_\circ+\ell_{L_1,T}+\ell_{L_2,T}]$:
\begin{equation}
\label{case2-bound}
\abs{\phi_T^\xi(k,k_\circ,\xi_\circ)} \leq \kappa_1(\abs{\xi_\circ}) \leq \kappa_1(\abs{\xi_\circ})e^{L_1(\abs{\xi_\circ})+L_2 \circ \kappa_1(\abs{\xi_\circ})}e^{-(k-k_\circ)T}  \ .
\end{equation}

Using definitions of $\ell_{L_i,T}$ and the fact that $T<T^*<1$, we obtain that
\begin{equation}
L_i-1 \leq T\ell_{L_i,T} \leq L_i, \ i=1,2 \ .
\end{equation}
Moreover, the definitions of $L_1,\nu$ we get that for all $k \geq k_\circ+\ell_{L_1,T}$,
\begin{equation}
\label{z-bound}
\abs{\phi_T^z(k,k_\circ,z_\circ)} \leq \max\{ \beta_z(\Delta,T\ell_{L_1,T}),\nu\} \leq \max\{ \beta_z(\Delta,L_1(\Delta)-1),\nu\} \leq c_1 \ .
\end{equation}
Furthermore, by letting $k_1 := k_\circ+\ell_{L_1,T}$, $k_2 := k_\circ+\ell_{L_1,T}+\ell_{L_2,T}$, $x_1 :=\phi_T^x(k_1,k_\circ,x_\circ,\omega^z_{[k_\circ,k_1)})$ and $x_2=\phi_T^x(k_2,k_\circ,x_\circ,\omega^z_{[k_\circ,k_2)})$ we obtain that 
\begin{eqnarray}
\label{x-bound}
\abs{\phi_T^x(k_2,k_\circ,x_\circ,\omega^z_{[k_\circ,k_2)})} & = & \abs{\phi_T^x(k_2,k_1,x_1,\omega^z_{[k_1,k_2)})} \nonumber \\
 & \leq & \max\{ \beta_x(\abs{x_1},(k_2-k_1)T), \nu_1 \}+\epsilon \nonumber \\
 & \leq & \max\{ \beta_x(\kappa_1(\Delta), T \ell_{L_2,T}),\nu_1\}+c_1/2 \nonumber \\
 & \leq & \max\{ \beta_x(\kappa_1(\Delta), L_2 \circ \kappa_1(\Delta)-1),\nu_1 \}+c_1/2 \nonumber \\
 & \leq & \max\{ c_1/2,c_1/2\}+c_1/2 \ = \ c_1 \ .
\end{eqnarray}
Hence, using (\ref{z-bound}), (\ref{x-bound}), Lemmas \ref{lem2} and \ref{lem-kl} and the construction of $\tilde{\beta}$ we have that for all $k \geq k_2$ :
\begin{eqnarray}
\abs{\phi_T^\xi(k,k_\circ,\xi_\circ)} & = &  \abs{\phi_T^\xi(k,k_2,\xi_2)} \nonumber \\
  & \leq & \max\{ \bar{\beta}(\abs{\xi_2},(k-k_2)T), \nu_1 \} \nonumber \\
  & \leq & \max\{ \bar{\beta}(\kappa_1(\abs{\xi_\circ}),(k-k_\circ-\ell_{L_1,T}-\ell_{L_2,T})T), \nu_1 \} \nonumber \\
 & \leq & \max\{ \bar{\beta}(\kappa_1(\abs{\xi_\circ}),(k-k_\circ)T-L_1-L_2), \nu_1 \} \nonumber \\
 & \leq & \max\{ \tilde{\beta}(\kappa_1(\abs{\xi_\circ}),(k-k_\circ)T), \nu \} \nonumber \ ,
\end{eqnarray}
which proves the result with the defined $\beta$.
\null \hfill \null $\blacksquare$

\subsection{Proof of Theorem \ref{th2}}
We again prove only semiglobal results since global results follow the same steps with minor changes.

\begin{lemma}
\label{lem:it:uos} 
Suppose that all conditions of Theorem \ref{th2} are satisfied. Let $V_T$ come from item 1 of the Theorem and define
\begin{equation}
\label{uos:y:def}
 y(k,\xi) :=  V_T(\,k+1, f_T(k,x,z) \, )  - V_T(\, k+1, f_{T}(k,x,0) \, )  \ .
\end{equation}
Then, there exists $\beta_y \in {\cal KL}$ such that for all $(\Delta,\nu)$ there exists $T^*>0$ such that for all $k_\circ \geq 0$, $\xi(k_\circ)=\xi_\circ$ with $\abs{\xi_\circ} \leq \Delta$ and $T \in (0,T^*)$ we have that 
\begin{equation}
\label{eq:fact:uos}
\abs{y(k,\phi_T^\xi(k,k_\circ,\xi_\circ))} \leq   T \cdot \max\{ \beta_y(\abs{\xi_\circ},(k-k_\circ)T) , \nu \} \,.
\end{equation}
for all $k \geq k_\circ$. \null \hfill \null $\square$
\end{lemma}

\noindent {\bf Proof. } 
Let $(\Delta,\nu)$ be given. Let $\gamma_1,\gamma_3 \in {\cal K}_\infty$, $\gamma_2 \in {\cal N}$ and $T_1^*>0$ come from Assumption \ref{ass2}. Let $\kappa$ and $c$ come from item 3 of the Theorem and let $\Delta$ generate $T_2^*>0$ via the same item. Let $L \in {\cal N}$ and $\Delta_1:=\gamma_1 (\kappa(\Delta)+c)$ generate $T_3^*>0$ via item 1 of the Theorem. Let $\beta_z \in {\cal KL}$ come from item 2 of the Theorem. Define
\begin{eqnarray}
\beta_y(s,t) & := & L \circ \gamma_1 (\kappa(s)+c) \cdot \gamma_2(\kappa(s)+c) \cdot \gamma_3 \circ \beta_z(s,t) \nonumber \\
\nu_1 & := & \gamma_3^{-1}\left( \frac{\nu}{L\circ\gamma_1(\kappa(\Delta)+c)\cdot\gamma_2(\kappa(\Delta)+c)} \right) \ . \nonumber 
\end{eqnarray}
Let $(\Delta,\nu_1)$ generate $T_4^*>0$ via the item 2 of Theorem. Let $T^*:=\min\{ T_1^*,T_2^*,T_3^*,T_4^* \}$ and consider arbitrary fixed $k_\circ \geq 0$, $\xi(k_\circ) = \xi_\circ$ with $\abs{\xi_\circ} \leq \Delta$ and $T \in (0,T^*)$. Then, using item the definition of $y$, Lipschitz condition on $V_T$ and Assumption \ref{ass2} we can write
\begin{eqnarray}
\abs{y(k,\xi)} & = & \abs{V_T(\,k+1, f_T(k,x,z) \, )  - V_T(\, k+1, f_{T}(k,x,0))} \nonumber \\
  & \leq & T \cdot L\left( \max\{ f_T(k,x,z), f_T(k,x,0) \} \right) \abs{f_T(k,x,z)-f_T(k,x,0)} \nonumber \\
  & \leq & T \cdot L \circ \gamma_1(\abs{\xi}) \cdot \gamma_2(\abs{x}) \cdot \gamma_3(\abs{z}) \ ,
\end{eqnarray}
for all $\abs{\xi} \leq \Delta_1$, $T \in (0,T^*)$. Finally, using items 2 and 3 of the Theorem and the definitions of $\beta_y$ and $\nu_1$ we can write
\begin{eqnarray}
\abs{y \left(k,\phi_T^\xi(k,k_\circ,\xi_\circ) \right) } & \leq & T \cdot L \circ \gamma_1 \left(\abs{\phi_T^\xi(k,k_\circ,\xi_\circ)} \right) \cdot \gamma_2 \left(\abs{\phi_T^\xi(k,k_\circ,\xi_\circ)} \right) \cdot \gamma_3(\abs{\phi_T^z(k,k_\circ,z_\circ)}) \nonumber \\
 & \leq & T \cdot L \circ \gamma_1(\kappa(\abs{\xi_\circ})+c) \cdot \gamma_2(\kappa(\abs{\xi_\circ})+c) \cdot \gamma_3(\max\{ \beta_z(\abs{z_\circ},(k-k_\circ)T),\nu_1\}) \nonumber \\
 & \leq & T \cdot \max\{ \beta_y(\abs{\xi_\circ},(k-k_\circ)T),\nu \} \ ,
\end{eqnarray}
for all $k \geq k_\circ$.

\noindent {\bf Proof of Theorem \ref{th2}:} 
We only show that there exists $\beta_x \in {\cal KL}$ and $\gamma \in {\cal K}$ such that for any $(\Delta,\nu)$ there exists $T^*>0$ such that for all $k_\circ \geq 0$, $\xi(k_\circ)=\xi_\circ$ with $\abs{\xi_\circ} \leq \Delta$ and $T \in (0,T^*)$ we have that:
\begin{equation} \label{x-conc}
\abs{\phi_T^x(k,k_\circ,x_\circ,\omega^z)} \leq \max\{ \beta_x(\abs{x_\circ},(k-k_\circ)T) , \nu \}+\gamma \left( \norm{\omega^z_{[k_\circ,k)}} \right) \ ,
\end{equation}
for all $k \geq k_\circ$, where $u(k) := \frac{y(k,\phi_T^\xi(k,k_\circ,\xi_\circ))}{T}$ and $\omega^z$ is a sequence of $u(k)$'s on the interval $[k_\circ,k)$. Then, the conclusion of the Theorem follows directly from Lemma \ref{lem-kl}, item 2 of the Theorem and Lemma \ref{lem:it:uos}. 

Let $\alpha_1,\alpha_2,\alpha_3$ come from item 1 of the Theorem. Define $\alpha(\cdot):=\frac{1}{2} \alpha_3 \circ \alpha_2^{-1}(\cdot)$ and consider the differential equation:
$$
\dot{\zeta} = -\alpha(\eta) \ , \zeta(0)=\zeta_\circ \ .
$$
Without loss of generality we can assume that $\alpha$ is locally a Lipschitz function (see \cite[footnote on pg. 139]{KHALIL}) and from \cite[Lemma 2.5]{SONWANSCL} we have that there exists $\beta \in {\cal KL}$ such that for all $0 \leq \zeta_\circ$ the solution of the differential equation exists, it is unique and satisfies:
$$
\zeta(t)=\beta(\zeta_\circ,t-t_\circ) \ , \ \forall \ t \geq t_\circ \geq 0.
$$
Define $\tilde{\alpha}(s):=\alpha^{-1}(s)+s$ and:
\begin{eqnarray}
\beta_y(s,t) & := & \alpha_1^{-1}\left( \beta(\alpha_2(s),t) \right) \nonumber \\
\gamma (s) & := & \alpha_1^{-1} \circ \tilde{\alpha}(2 s) \ .
\end{eqnarray}
Let $(\Delta,\nu)$ be given. Let $\Delta_1:=\kappa(\Delta)+c$, where $\kappa,c$ come from the item 3 of Theorem. Let $\nu_1:=\frac{1}{2} \tilde{\alpha}^{-1} \circ \alpha_1(\nu)$. Let $(\Delta_1,\nu_1)$ generate $T^* \in (0,1)$ be so that the items 1 and 3 of Theorem hold. Consider arbitrary fixed $k_\circ \geq 0$, $\xi(k_\circ) = \xi_\circ$ with $\abs{\xi_\circ} \leq \Delta$ and $T \in (0,T^*)$. Note that the item 3 of Theorem guarantees that the solutions of (\ref{sys1}), (\ref{sys2}) exist for all $k \geq k_\circ \geq 0$ and satisfy $\abs{\phi_T^\xi(k,k_\circ,\xi_\circ)} \leq \Delta_1$. We can write
\begin{eqnarray}
V_T(k+1, f_T(k,x,z)\,) - V_T(k,x) & = &  V_T(\,k+1, f_T(k,x,0) \,) - V_T(k,x) + \nonumber \\
 & & V_T(\,k+1, f_T(k,x,z) \, ) -  V_T(k+1, f_T(k,x,0) \,) \nonumber \\
 & = & V_T(\,k+1, f_T(k,x,0) \,) - V_T(k,x) + Tu(k) \ .
\label{seven}
\end{eqnarray}
Then, defining $v(k):=V_T(k,\phi_T^x(k,k_\circ,x_\circ,\omega^z_{[k_\circ,k)}))$, using (\ref{newa1a:VT:1}), (\ref{newa1a:VT:2}), the following holds along trajectories
\begin{eqnarray}
\nonumber
v(k+1) - v(k) & \leq & -T 2 \alpha(v(k)) +T [\nu_1 + u(k)] 
\label{eight}   \,
\end{eqnarray}
for all $k \geq k_\circ \geq 0$. From this inequality we can see that the following hold true for any $k\geq k_\circ\geq 0$,
\begin{eqnarray}
\label{iss}
&\displaystyle  v(k) \geq \alpha^{-1} \left(\nu_1+\norm{\omega^z_{[k_\circ,k)}} \right)  \ \Longrightarrow\ v(k+1)- v(k) \leq -T \alpha(v(k)) & \\[1mm]
\label{eq:growth}
&\displaystyle v(k+1) \leq v(k) + T \left(\nu_1+\norm{\omega^z_{[k_\circ,k)}} \right) \,.&
\end{eqnarray}
Define next $\zeta(t):= v(k) + \left(\displaystyle\frac{t}{T} -k \right) ( v(k+1) - v(k) )$ for all $t \in [kT, (k+1)T\,)$, $k\geq k_\circ\geq 0$. Notice that since $\zeta(t)$ is a linear interpolation from $v(k)$ to $v(k+1)$, which are always non-negative, we have that
\begin{equation}
\label{bndsonzeta}
0 \leq \zeta(t) \leq \max\{v(k),\,v(k+1)\} \qquad \mbox{for any  \ \  $t \in [kT, (k+1)T\,)$, $k\geq k_\circ\geq 0$}\,.
\end{equation}
Moreover, since trajectories of the system (\ref{sys1}), (\ref{sys2}) are defined for all $k \geq k_\circ \geq 0$, the variable $\zeta(t)$ is defined for all $t \geq t_\circ \geq 0$. Introduce $\bar{u}(t):=u(k), \forall t \in [kT,(k+1)T), k \geq k_\circ \geq 0$. Using definition of $\tilde{\alpha}$, $T^*<1$, (\ref{iss}), (\ref{eq:growth}) and (\ref{bndsonzeta}) we can write for all $t \in [kT,(k+1)T), k \geq k_\circ \geq 0$:
\begin{eqnarray}
\zeta(t) \geq \tilde{\alpha}(\nu_1+\norm{\bar{u}_{[t_\circ,t)}}) & \Longrightarrow & v(k) \geq \alpha^{-1} \left(\nu_1+\norm{\omega^z_{[k_\circ,k)}} \right) \nonumber \\
  & \Longrightarrow & v(k+1) - v(k)  \leq  -T \alpha(v(k)) \nonumber \\
   & \Longrightarrow & \dot{\zeta}(t) \leq -T \alpha(v(k)) \nonumber \\
   & \Longrightarrow & \dot{\zeta}(t) \leq -T \alpha(\zeta(t)) \ ,
\end{eqnarray}
where the last inequality follows from the fact that $v(k+1) \leq v(k)$ and hence $\zeta(t) \leq v(k)$ for all $t \in [kT,(k+1)T), k \geq k_\circ \geq 0$. Using \cite[Corollary 5.1 and 5.2]{KHALIL}, we have that:
$$
\zeta(t) \leq \max\{ \beta(\zeta_\circ,t-t_\circ), \tilde{\alpha} (\nu_1+\norm{\bar{u}_{[t_\circ,t)}}) \} \ ,
$$
and using $\phi_T^x(k,k_\circ,x_\circ,\omega^z_{[k_\circ,k)}) = \zeta(kT)$, $\bar{u}(kT)=u(k)$, $\beta_y$ and $\nu_1$, (\ref{newa1a:VT:1}) and the fact that for arbitrary $\alpha \in {\cal K}_\infty$ we have $\alpha(r+s) \leq \alpha(2s)+\alpha(2r), \forall r,s \geq 0$ we can write
\begin{eqnarray}
\abs{\phi_T^x(k,k_\circ,x_\circ,\omega^z_{[k_\circ,k)})} & \leq & \alpha^{-1}_1 \left( \max \left\{ \beta(\alpha_2(\abs{x_\circ}),(k-k_\circ)T), \tilde{\alpha} (\nu_1+\norm{\omega^z_{[k_\circ,k)}}) \right\} \right) \nonumber \\
 & \leq & \max\left\{ \alpha^{-1}_1(\beta(\alpha_2(\abs{x_\circ}),(k-k_\circ)T)),  \alpha_1^{-1} \circ \tilde{\alpha} (2\nu_1) +\alpha_1^{-1} \circ \tilde{\alpha}(\norm{\omega^z_{[k_\circ,k)}}) \right\} \nonumber \\
 & \leq & \max\{ \beta_y(\abs{x_\circ},(k-k_\circ)T),\nu\}+\gamma(\norm{\omega^z_{[k_\circ,k)}}) \ ,
\end{eqnarray}
for all $k \geq k_\circ \geq 0$.

\section{Conclusions}
\label{sec:concl}
We have established necessary and sufficient conditions for semiglobal practical asymptotic stability of time-varying parameterized cascades. Parameterized systems arise naturally when an approximate discrete-time model of a sampled-data plant is used for controller design. Our results provide the controller designer with a range of tools that can be used for a systematic digital controller design based on approximate discrete-time models. The utility of our results was illustrated with a case study, where we have obtained a controller that performs better than the emulated (discretized) continuous-time controller designed in \cite{MOBCAR} and implemented using a sampler and zero-order-hold. 


\bibliographystyle{plain}
\bibliography{refs}


\appendix

\section{Some technical proofs}

\subsubsection*{Proof of Proposition \ref{lem3}:} 
The result is only proved for USC since the proof of UGC follows the same steps. Let $\beta_x \in {\cal KL}$ come from SP-UAS. Let $(\Delta,\nu)$ be given and let $T_1^*>0$ be generated using SP-UAS. Without loss of generality assume that $\Delta>\nu$. Let $\Delta_1:=\beta_x(\Delta,0)+1$. Let $\Delta_2>0$ be arbitrary. Let $(\Delta_1,\Delta_2)$ generate $K,T_2^*>0$ using conditions of the Lemma. Let $T^*:=\min \{ T_1^*,T_2^* \}$. Consider arbitrary numbers $\eta \in (0,\Delta)$, $\epsilon,L>0$ and define
$$
\mu:= \min \left\{ \frac{\min\{ \epsilon,1\} }{e^{KL}-1},\Delta_2 \right\} \ .
$$
Consider arbitrary $\abs{x_\circ} \leq \eta$, $\norm{\omega^z} \leq \mu$, $T \in (0,T^*)$. Note first that without loss of generality we may assume that $\beta_x$ satisfies \rref{beta0}. Consequently, SP-UAS implies that whenever $\abs{x_\circ} \leq \Delta$, $T \in (0,T^*)$ then $\abs{\phi_T^x(k,k_\circ,x_\circ,0)} \leq \max\{\beta_x(\Delta,0),\nu\} \leq \beta_x(\Delta,0) <\Delta_1$ for all $k \geq k_\circ \geq 0$. 

With the goal of showing contradiction suppose that there exists $k_1^* \in (k_\circ,k_\circ+\ell_{L,T})$ such that $\abs{\phi_T^x(k,k_\circ,x_\circ,\omega^z_{[k_\circ,k)})} \leq \Delta_1$ for all $k \in [k_\circ,k_1^*)$ and $\abs{\phi_T^x(k_1^*,k_\circ,x_\circ,\omega^z_{[k_\circ,k_1^*)})}>\Delta_1$. Then, we can write
\begin{eqnarray}
\label{lemma1:last}
\abs{\phi_T^x(k,k_\circ,x_\circ,\omega^z_{[k_\circ,k)})-\phi^x_T(k,k_\circ,x_\circ,0)} & \leq & KT \sum_{j=0}^{k-k_\circ-1} (1+KT)^j \norm{\omega^z_{[k_\circ,k)}} \nonumber \\
 & \leq & [(1+KT)^{k-k_\circ}-1] \norm{\omega^z} \nonumber \\
 & \leq & [e^{KT(k-k_\circ)}-1]\norm{\omega^z} \nonumber \\
 & \leq & [e^{KL}-1] \norm{\omega^z} \nonumber \\
 & \leq & \min\{ \epsilon ,1 \} 
\end{eqnarray}
for all $k \in \left[k_\circ,k_1^* \right]$. This contradicts the assumption that $\abs{\phi_T^x(k_1^*,k_\circ,x_\circ,\omega^z_{[k_\circ,k_1^*)})}>\Delta_1$ since
\begin{eqnarray*}
\abs{\phi_T^x(k_1^*,k_\circ,x_\circ,\omega^z_{[k_\circ,k_1^*)})} & \leq & \abs{\phi^x_T(k_1^*,k_\circ,x_\circ,0)} +\abs{\phi_T^x(k_1^*,k_\circ,x_\circ,\omega^z_{[k_\circ,k_1^*)})-\phi^x_T(k_1^*,k_\circ,x_\circ,0)} \\
& \leq & \beta_x(\Delta,0)+1=\Delta_1 \ . 
\end{eqnarray*}
Hence  (\ref{lemma1:last}) holds for all $k \geq k_\circ+\ell_{L,T}$, which completes the proof.
\null \hfill \null $\blacksquare$

\subsubsection*{Proof of  Claim 1:}
Let $T_1^*>0$ come from inequalities (\ref{c1c2V}) and (\ref{93}). Let $\mu,L,T_2^*$ come from the PE condition. Let $w_M$ come from the conditions of the proposition. Let $T_3^*>0$  and $T_4^*>0$ be such that 
\begin{equation}\label{T*1}
\frac{T}{1-e^{-T}} \leq 2, \qquad \forall \, T \in (0,T^*_3); \qquad \frac{1-e^{-T}}{T} \geq \frac{1}{2}, \qquad \forall\, T \in (0,T^*_4) \ .
\end{equation}
Let $T_5^*:= \frac{\mu e^{-L}}{4(1-e^{-L})}$ and, finally, define $T^*:=\min\{T_1^*,T_2^*,T_3^*,T_4^*,T_5^*\}$. Consider arbitrary $x \in \reals^2$, $k \geq 0$ and $T \in (0,T^*)$. The first inequality in (\ref{w1}) holds using our choice of $T_3^*$:
\begin{eqnarray}
T \sum_{i=k}^{\infty} e^{(k-i)T}\omega_{r_i}^2 & \leq & w_M^2 T \sum_{j=0}^{\infty} \left( e^{-T} \right)^j \ = \ w_M^2 \frac{T}{1-e^{-T}} \ \leq \  2 w_M^2 \ =: c_3 \ . \label{c3}
\end{eqnarray}
The second inequality in (\ref{w1})  follows from the PE condition. In particular, denote $k_j:=k+j\ell_{L,T}$ and then we can write:
\begin{eqnarray}
T\sum_{i=k}^{\infty} e^{(k-i)T}\omega_{r_i}^2 & = & \sum_{j=0}^{\infty} T \sum_{i=k+j\ell_{L,T}}^{k+(j+1)\ell_{L,T}} e^{(k-i)T}\omega_{r_i}^2 \nonumber \\
 & \geq & \sum_{j=0}^{\infty} e^{-jT\ell_{L,T}} \left( T \sum_{i=k_j}^{k_j+\ell_{L,T}} e^{(k_j-i)T}\omega_{r_i}^2 \right) \nonumber \\
 & \geq & \sum_{j=0}^{\infty} e^{-jL} e^{-\ell_{L,T}T} \left( T \sum_{i=k_j}^{k_j+\ell_{L,T}} \omega_{r_i}^2 \right)  \nonumber \\
 & \geq & \frac{e^{-L}}{1-e^{-L}} \mu \ =: \ c_4 \ . \label{c4}
\end{eqnarray}
We show next that (\ref{w2}) holds. Using similar calculations as in (\ref{c3}), (\ref{c4}) and our choice of $T_4^*$ and $T_5^*=\frac{c_4}{4}$ we have:
\begin{eqnarray*}
\frac{W_T(k+1,y_e-\omega_{r_k}x_e)-W_T(k,y_e)}{T} & = & - \sum_{i=k+1}^{\infty} e^{(k+1-i)T}\omega_{r_i}^2 (y_e-T\omega_{r_k}x_e)^2 + \sum_{i=k}^{\infty} e^{(k-i)T}\omega_{r_i}^2 y_e^2 \nonumber \\
 & = & \omega_{r_k}^2 y_e^2 - \sum_{i=k+1}^{\infty} e^{(k+1-i)T}\omega_{r_i}^2 \left[ (y_e-T\omega_{r_k}x_e)^2 -y_e^2\right] \nonumber \\
 &  & - (1-e^{-T}) \sum_{i=k+1}^{\infty} e^{(k+1-i)T}\omega_{r_i}^2 y_e^2 \nonumber \\
 & \leq & \omega_{r_k}^2 y_e^2 + \left( \omega^2_M T \sum_{j=0}^\infty e^{-j T}\right) \left[ 2 w_M \abs{x_e}\abs{y_e} + T w_M^2 x_e^2 \right] \nonumber \\
 & & - \frac{1-e^{-T}}{T} T \sum_{i=k+1}^{\infty} e^{(k+1-i)T}\omega_{r_i}^2 y_e^2 \nonumber \\
 & \leq & \omega_{r_k}^2 y_e^2 + \frac{w_M^2 T}{1-e^{-T}} (2w_M \abs{x_e} \abs{y_e}+T w_M^2 x_e^2) - \frac{c_4}{2} y_e^2  \ . \label{dw2}
\end{eqnarray*}
Finally, using {\rref{T*1} and completing the} squares {we obtain that} there exist $K_2,\tilde \alpha_y>0$ such that 
\begin{equation}\label{poclednie!}
\frac{W_T(k+1,y_e-\omega_{r_k}x_e)-W_T(k,y_e)}{T}  \leq  \omega_{r_k}^2 y_e^2+K_2 x_e^2 - \tilde \alpha_y y_e^2 \ ,
\end{equation}
which completes the proof of the claim.
\null\hfill\null $\blacksquare$

\newpage

\subsection{Proof of Proposition \ref{prop:iisns-lyap-new}}
\label{sec:proofprop2}

First, note that using (\ref{iisns-new-eq2}), (\ref{iisns-eq3}) we can write
\begin{eqnarray}
\label{eqa}
V_T(k+1,f_T(k,x,z))-V_T(k,x) & = & V_T(k+1,f_T(k,x,0))-V_T(k,x) \nonumber \\
 & & +V_T(k+1,f_T(k,x,z))-V_T(k+1,f_T(k,x,0)) \nonumber \\
 & \leq &  T \tilde{\gamma}_1(\abs{z}) \varphi(V_T(k,x)) + T \tilde{\gamma}_2(\abs{z}) . \nonumber 
\end{eqnarray}
Assume for the time-being that from this inequality and the conditions of the Proposition, there exist $\alpha_1,\alpha_2 \in {\cal K}_\infty$, $c>0$ and for each $T \in (0,T^*)$ there exists $W_T: \reals_{\geq 0} \times \reals^{n_x} \rightarrow \reals_{\geq 0}$ such that all $x \in \reals^n$, $k \geq 0$ and $T \in (0,T^*)$ we have that 
\begin{eqnarray}
\alpha_1(\abs{x}) \ \leq \ W_T(k,x)  \ & \leq & \alpha_2(\abs{x})+c \label{eq-ufc1-a} \\
W_T(k+1,f_T(k,x,z))-W_T(k,x) & \leq & T \mu(\abs{z})   \label{eq-ufc2-a} \,.
\end{eqnarray} 
Then, defining $w(k):=W_T(k,\phi_T^x(k,k_\circ,x_\circ,\omega^z_{[k_\circ,k)}))$ and using (\ref{eq-ufc2-a}) we can write
$$
w(k+1) \leq w(k)+ T \mu(\abs{\omega^z_{[k_\circ,k)}}) , \qquad \forall k \geq k_\circ \geq 0 \ 
$$
and consequently, by induction it follows that
$$
w(k) \leq w(k_\circ)+T \sum_{i=k_\circ}^{k-1} \mu(\abs{\omega^z_{[k_\circ,k)}}) , \qquad \forall k \geq k_\circ \geq 0 \ .
$$
Next, using (\ref{eq-ufc1-a}) we obtain that 
\begin{eqnarray}
\alpha_1\left(\abs{\phi_T^x(k,k_\circ,x_\circ,\omega^z_{[k_\circ,k)})}\right) & \leq & \alpha_2(\abs{x_\circ})+c+T\sum_{t=k_\circ}^{k-1} \mu(\abs{\omega^z_{[k_\circ,k)}}) , \qquad \forall k \geq k_\circ \geq 0 \ ,
\end{eqnarray}
hence, UGB follows from \rref{eq:summability} and the fact that $\alpha_1,\,\alpha_2\in\cKinfty$.
 
It is only left to prove that under the conditions of the proposition \rref{eq-ufc1-a} and \rref{eq-ufc2-a} hold. To that end let $\tilde{\alpha}_1,\tilde{\alpha}_2,\varphi,\tilde{\gamma}_1,\tilde{\gamma}_2$ come from Assumption \ref{mars}. Define
$$
q(s) : = \left\{ 
\begin{array}{ll}
\frac{1}{\varphi(1)} , & s \leq 1 \\
\frac{1}{\varphi(s)} , & s > 1
\end{array}
\right. \ .
$$
Let $\rho \in {\cal K}_\infty$ be defined as 
$$
\rho(s) := \int_0^s q(\tau)d\tau \ .
$$
Note that $\rho \in {\cal K}_\infty$ and $\rho$ is differentiable with $\frac{d\rho}{ds}(s)=q(s)$. Also, $q(\cdot)$ is non-increasing since $\varphi \in {\cal K}_\infty$. Introduce $W_T(k,x):=\rho(V_T(k,x))$. Then, it is clear that (\ref{eq-ufc1-a}) holds with $\alpha_1(s):=\rho \circ \tilde{\alpha}_1(s)$, $\alpha_2(s)=\rho \circ 2 \tilde{\alpha}_2(s)$ and $c=\rho(2\tilde{c})$. We now show that (\ref{eq-ufc2-a}) holds.

We consider two cases. If $V_T(k+1,f_T(k,x,z)) \leq V_T(k,x)$, then since $\rho \in {\cal K}_\infty$, we have:
\begin{equation}
\label{w-new}
W_T(k+1,F_T(k,x,z))-W_T(k,x) \leq 0  \leq T \left( \tilde{\gamma}_1(\abs{z}) +\tilde{\gamma}_2(\abs{z}) \frac{1}{\varphi(1)} \right) \ .
\end{equation}
On the other hand, if $V_T(k+1,f_T(k,x,z)) > V_T(k,x)$, then using the Mean Value Theorem we can write
\begin{eqnarray}
W_T(k+1,f_T(k,x,z))-W_T(k,x) & = & q(V^*) [V_T(k+1,f_T(k,x,z)) - V_T(k,x)] \nonumber \\
 & \leq & T q(V^*)[\tilde{\gamma}_1(\abs{z})\varphi(V_T(k,x))+\tilde{\gamma}_2(\abs{z})]
\end{eqnarray}
where $V_T(k,x) <V^* <V_T(k+1,f_T(k,x,z))$. Since $q$ is a non increasing function, we have
$$
q(V_T(k,x))) \geq q(V^*) \ ,
$$
and hence
\begin{eqnarray}
W_T(k+1,f_T(k,x,z))-W_T(k,x) & \leq & T q(V_T(k,x))[\tilde{\gamma}_1(\abs{z})\varphi(V_T(k,x))+\tilde{\gamma}_2(\abs{z})] \ . \nonumber 
\end{eqnarray}
Now, if $V_T(k,x) \leq 1$ then we have that
\begin{eqnarray}
\label{w-new1}
W_T(k+1,f_T(k,x,z))-W_T(k,x) & \leq & T \frac{\tilde{\gamma}_1(\abs{z})\varphi(V_T(k,x))+\tilde{\gamma}_2(\abs{z}) }{\varphi(1)} \nonumber \\
 & \leq & T \frac{\tilde{\gamma}_1(\abs{z}) \varphi(1)+\tilde{\gamma}_2(\abs{z})}{\varphi(1)} \nonumber \\
 & = & T \left( \tilde{\gamma}_1(\abs{z}) +\tilde{\gamma}_2(\abs{z}) \frac{1}{\varphi(1)} \right)  \ .
\end{eqnarray}
Otherwise, if $V_T(k,x) \geq 1$ then, 
\begin{eqnarray}
\label{w-new2}
W_T(k+1,f_T(k,x,z))-W_T(k,x) & \leq & T \frac{\tilde{\gamma}_1(\abs{z})\varphi(V_T(k,x))+\tilde{\gamma}_2(\abs{z})}{\varphi(V_T(k,x))} \nonumber \\
 & \leq & T\left(\tilde{\gamma}_1(\abs{z})+\frac{\tilde{\gamma}_2(\abs{z})}{\varphi(V_T(k,x))} \right)  \nonumber \\
 & = & T \left( \tilde{\gamma}_1(\abs{z}) +\tilde{\gamma}_2(\abs{z}) \frac{1}{\varphi(1)} \right)  \ .
\end{eqnarray}
The proof is completed by combining (\ref{w-new}), (\ref{w-new1}) and (\ref{w-new2}).
\null \hfill \null $\blacksquare$ 

\end{document}